\documentclass[12pt]{amsart}
\usepackage{amscd,amsmath,amssymb,amsfonts}
\theoremstyle{plain}
\newtheorem{thm}{Theorem}
\newtheorem{lem}[thm]{Lemma}
\newtheorem{cor}[thm]{Corollary}
\newtheorem{prop}[thm]{Proposition}

\theoremstyle{definition}
\newtheorem{defn}[thm]{Definition}
\newtheorem{remark}[thm]{Remark}

\newtheorem{notation}[thm]{Notation}
\numberwithin{thm}{section}
\numberwithin{equation}{section}

\newcommand{\e}{\epsilon}

\newcommand{\Sym}{{\rm Sym}}

\newcommand{\id}{{\rm id}}

\newcommand{\Tr}{{\rm Tr}}

\newcommand{\ov}{\widetilde}

\newcommand{\sA}{{\mathcal A}}

\newcommand{\sC}{{\mathcal C}}
\newcommand{\sD}{{\mathcal D}}
\newcommand{\sE}{{\mathcal E}}
\newcommand{\sF}{{\mathcal F}}
\newcommand{\sG}{{\mathcal G}}
\newcommand{\sH}{{\mathcal H}}

\newcommand{\sK}{{\mathcal K}}
\newcommand{\sL}{{\mathcal L}}
\newcommand{\sM}{{\mathcal M}}
\newcommand{\sN}{{\mathcal N}}
\newcommand{\sO}{{\mathcal O}}

\newcommand{\sR}{{\mathcal R}}

\newcommand{\sU}{{\mathcal U}}
\newcommand{\sV}{{\mathcal V}}
\newcommand{\sW}{{\mathcal W}}

\newcommand{\E}{{\mathbb E}}

\newcommand{\M}{{\mathbb M}}

\newcommand{\Z}{{\mathbb Z}}

\begin{document}

\title
[Hitchin's connection and 2nd Order Operators]
{Hitchin's connection and Differential Operators
with Values in the Determinant Bundle}
\author{Xiaotao Sun}
\address{Institute of Mathematics, AMSS, Chinese Academy of Sciences, China}
\email{xsun@mail2.math.ac.cn}
\address{Department of Mathematics, The University of Hong Kong, Hong Kong}
\email{xsun@maths.hku.hk}
\author{I-Hsun Tsai}
\address{Department of Mathematics, National Taiwan University, Taipei, 
Taiwan}
\email{ihtsai@math.ntu.edu.tw}
\date{}
\subjclass{}
\thanks{Sun is supported by the outstanding young grant of NFSC at number
10025103 and a grant of RGC of Hong Kong at RGC code 7131/02P.  
Tsai has been supported in part by 
a grant from National Science Council of Taiwan.}     
\maketitle

\section{Introduction}

For any smooth family $\pi: X \to S$ of curves 
and a vector bundle $E$ on $X$,  
A. Beilinson and V. Schechtman 
defined in \cite{BS} (see also \cite{ET}) a so-called trace complex
$\,^{tr}\sA_E^{\bullet}$ on $X$, together with an algebra structure on it, 
such that
$R^0\pi_*(\,^{tr}\sA_E^{\bullet})$ is canonically isomorphic to
the Atiyah algebra $\sA_{\lambda_E}$ of the determinant
bundle $\lambda_E=\det R\pi_*E$ of the family 
(See Proposition \ref{prop4.3} for details). It is generalized by
Y.-L. Tong and the second author (\cite{TT}) 
to the families $\pi:\ov X\to\ov S$ of stable curves. 
On the other hand, in his fundamental article \cite{H},
N. Hitchin constructed a projective connection on
the relative sections of the determinant bundle on the moduli of
bundles with trivial determinant on curves of fixed genus. 
In this paper, we present 
{\it a construction of Hitchin's connection and its
logarithmic extension 
within the framework of \cite{BS} and \cite{TT}}.  
Note that studying the behaviour of spaces of
generalized theta functions under degeneration of the curves was
already suggested in \cite{BK} and \cite{H}.

Let $\sM$ be the moduli stack of smooth curves of genus $g\ge 2$ and
$\sC\to\sM$ be the universal curve. Let $f: S\to \sM$ be the family
of moduli spaces of stable bundles of rank $r\ge 2$ with a fixed
determinant. Let $\pi:X=\sC\times_{\sM}S\to S$, $E$ a universal
bundle on $X$ and $\E:=\sE nd^0(E)$. In view of identification of
tangent bundle $T_{S/\sM}$ with $R^1\pi_*(\E)$, it
seems natural to ask whether there exists an
analogous theory for the side of differential
operators, namely for an identification
of sheaves of differential operators (with values
in the determinant bundle) with certain cohomology groups. 
The trace complex of \cite{BS} provides
the identification, but we found unfortunately it hard to work with
the trace complex in a direct way. To tackle this difficulty we find that
there exist {\it canonical
(locally free) sheaves} $\sG_E$, $S(\sG_E)$ whose relative cohomology 
coincides with that of the trace complex (at least in 
the moduli situation).  Thus we have canonical
identifications
$$\phi:R^1\pi_*(\sG_E)\cong \sD^{\le 1}_{S/\sM}(\lambda_{\E}),\quad
\phi:R^1\pi_*S(\sG_E)\cong \sD^{\le 1}_S(\lambda_{\E})$$
(cf. Theorem \ref{prop0.2}). This part is of independent interest in itself. 

Let $\Theta$ be the theta line bundle on $S$ and, for any $k$, 
write $\Theta^k=\lambda_{\E}^{\mu}$ for some $\mu\in\Bbb Q$. 
Using the above identifications
and taking cohomology for some short exact sequences, we get the commutative
diagram
$$\CD
\sD^{\le 1}_{S/\sM}(\Theta^k)@>\bullet\mu>>\sD^{\le 1}_{S/\sM}(\Theta^k)\\
@VVV               @VVV \\
\sD^{\le 1}_S(\Theta^k)@>\tilde\delta_{\rm H}>>
S^2\sD^{\le 1}_{S/\sM}(\Theta^k)\\
@V\sigma VV                       @V S^2\e_1VV\\
f^*T_{\sM} @>\delta_{\rm H}>>S^2T_{S/\sM}
\endCD$$
where the two vertical sequences are short exact sequences 
and $\delta_{\rm H}$, $\tilde\delta_{\rm H}$ are the connecting maps.
$\delta_{\rm H}$ is in fact
an {\it alternative of Hitchin's symbol map} (cf. \cite{R}, Section 4).

As nicely explained in \cite{GJ}, Hitchin's construction of
the connection on the vector bundle $\sV_k=f_*(\Theta^k)$
consists in defining, for any 
$v\in T_{\sM}(U)=f^*T_{\sM}(f^{-1}(U))$, a {\it heat operator} 
${\rm H}(v)\in H^0(f^{-1}(U),\sD^{\le 2}_S(\Theta^k))$  
such that its symbolic part is 
$\delta_{\rm H}(v)$. With help of above diagram, we are able to compare
the obstruction classes to lift $\delta_{\rm H}(v)$
to $\sD^{\le 2}_{S/\sM}(\Theta^k)$ and
to $S^2\sD^{\le 1}_{S/\sM}(\Theta^k)$ (cf. Proposion \ref{prop0.3}). 
Then we get a lifting  
${\rm H}(v)$. 
The uniqueness of ${\rm H}(v)$
follows from the vanishing of $f_*T_{S/\sM}$. It is important that we
do not use the vanishing of $R^1f_*\sO_S$, so that 
our method can be extended to some stable curves, where such a
vanishing may not hold. 

Let $\ov\sM\subset \ov{\sM_g}$
be the open set consisting of irreducible stable
curves and reducible curves with only one node.
Let $\ov\sC\to \ov\sM$ be the universal curve (extended
$\sC\to\sM$), and $\ov f: \ov S\to \ov\sM$ be the family of moduli spaces
of stable bundles with fixed determinant (cf. Notation \ref{prop2}).
Then, 
in the context of {\it log-geometry}
the similar formulation for $\pi:\ov X:=\ov\sC\times_{\ov\sM}\ov S\to\ov S$, 
via an appropriate generalization of the trace complex ( cf. \cite{TT}), 
can be achieved equally well.  
Thus another result of this paper is that there is
a {\it logarithmic projective connection} on
${\ov f}_*(\Theta^k)$. 
We remark that $R^1{\ov f}_*\sO_{\ov S}=0$ may not hold, but
$\ov{f}_*T_{\ov S/\ov\sM}=0$ remains true since $S$ is dense
in $\ov S$ (cf. Notation \ref{prop2}). 
The family $\ov f: \ov S\to \ov\sM$ is not proper. However,
we can prove that ${\ov f}_*(\Theta^k)$ is a {\it coherent} sheaf
on $\ov\sM$ (cf. Theorem \ref{prop9.2}).   We believe that 
the coherent property remains true on  
$\ov M_g$ of all stable curves; this is not, however, proved in the 
present paper.     

Basically our approach works equally well for
both of the general case \cite{H} ($g\ge 3$) and special cases
(e.g. $g=2$ and $r=2$) \cite{GJ}, with the difference 
arising from the fact that in $g=r=2$ case our approach 
shall work with a Quot scheme 
$\sR^{ss}$ which has a good quotient $\sR^{ss}\to S$ where $S$ is the   
family of moduli spaces of S-equivalence classes of {\it semistable}  
vector bundles.   As this will introduce additional details, 
we leave the case $g=r=2$ to future discussions.      

In Section 2, we construct Hitchin's connection and its logarithmic
extension under assumption of Theorem \ref{prop0.2}, 
Lemma \ref{prop0} and Theorem \ref{prop0.6}.
It is enough to read this section
for a reader who is only interested in Hitchin's connection.
Section 3 is devoted to the proof of Theorem \ref{prop0.2} and 
Lemma \ref{prop0}. In Section 4, we indicate the modifications of
arguments in Section 3 to prove Theorem \ref{prop0.6}. Most of
Section 4 is devoted to prove the {\it coherence} of $\ov{f}_*(\Theta^k)$
on $\ov\sM$. 

Remark that some ideas closely related 
to Section 3 of this paper were briefly discussed in 
Sections 9, 10 of \cite{G} where the subject matter 
is considered for moduli space of (stable) $G$-bundles, with 
$G$ being complex semisimple, connected and simply-connected.       
   
{\it Acknowledgements}.  We would like to express our deep 
gratitude to H\'el\`ene Esnault for the great impact of 
her ideas on this work.  Indeed numerous discussions with           
her in an earlier stage of this work, 
either directly allow us to grasp several key ideas, 
or contribute in an essential way to shaping 
ourselves towards the present approach.   
This joint work was initiated while both of us    
were visiting Essen (though a few topics involved
in this work had already been studied by three of us 
independently).  We wish to thank heartily  
H\'el\`ene Esnault and Eckart Viehweg 
for their warm hospitality and 
the excellent working conditions of Essen, that 
made our visit and this work possible.

\section {Heat operators and cohomology of sheaves}

As in Introduction, $\sM$ (resp. $\ov{\sM_g}$) denotes the moduli stack
of smooth (resp. stable) curves of genus $g\ge 2$,
and $\sC\to\sM$
(resp. $\ov\sC\to\ov{\sM_g}$) denotes the universal curve respectively.
By {\it moduli stack} in this article, we mean that we are working
on the {\it fine moduli spaces}. In particular, we assume that both
the moduli spaces $\sM$, $\ov{\sM_g}$ and the universal curves $\sC$,
$\ov\sC$ are smooth.  Fix a line bundle $\sN$ of relative degree
$d$ on $\ov\sC\to\ov{\sM_g}$, let $f:S\to\sM$ be the family of moduli
spaces of stable bundles of rank $r$ with fixed determinant
$\sN_b:=\sN|_{\sC_b}$ ($b\in\sM$) on $\sC_b$.
Let $\pi:X=\sC\times_{\sM}S\to S$ be the pull-back of
$\sC\to\sM$ via $f:S\to\sM$.

\begin{remark} If $(d, r)=1$,
there exists a universal bundle $E$ on $X$.
In general $S$ is a good quotient of a Hilbert quotient scheme,
denoted by $\sR^s$,
such that there is a universal bundle $E$ on
$X_{\sR^s}:=\sC\times_{\sM} \sR^s$.
$E$ may not descend to $X$, but objects such as
$\sE nd^0(E)$, $\sG_E$, $^{tr}\sA_E^{\bullet}$
(and other relevant constructions
that will be discussed later in this section) do descend.
Recall that a sheaf $\sF$ on $X_{\sR^s}$ descends
to $X$ if the action, of scalar automorphisms of $E$
(relative to $\sR^s$), on $\sF$ is trivial, e.g. \cite{HL}.
\end{remark}

{\it Without the danger of confusion
we will henceforth be working as if
a universal bundle $E$ exists on $X$.}

Let $\Theta$ be the theta line bundle on $S$.
For a sheaf $\sF$ on $X$,
$$\lambda_{\sF}=\bigotimes_{q\ge 0}{\rm det}\,R^q\pi_*\sF^{(-1)^q}
$$
denotes the Knudsen-Mumford determinant bundle on $S$.
As usual $T_S$, $K_S$
denote tangent bundle, canonical bundle;
$T_{S/\sM}$, $K_{S/\sM}$
denote the relative counterparts. Assume that $K_{S/\sM}=\Theta^{-\lambda}$.
Let $\sD^{\le i}_{S}(\Theta^k)$ be the sheaf of differential
operators on $\Theta^k$ of order $\le i$,
and by $\epsilon$ the symbol map of differential operators.
Let $\sW(\Theta^k):=\sD^{\le 1}_S(\Theta^k)+
\sD^{\le 2}_{S/\sM}(\Theta^k),$ one has an exact
sequence
\begin{gather}\label{eqsigma}
0\to \sD^{\le 2}_{S/\sM}(\Theta^k) \to
\sW(\Theta^k)
\xrightarrow{\sigma} f^*T_\sM \to 0.
\end{gather}

\begin{defn}  (cf. \cite{GJ}, 2.3.2.)
A {\it heat operator} ${\rm H}$
on $\Theta^k$ is an $\sO_S$-map
$f^* T_\sM\xrightarrow{\rm H} \sW(\sF) $
which, while composed with $\sigma$ above,
is the identity.  A {\it projective heat operator}
is ${\rm H}:T_{\sM}\to
f_*\sW(\Theta^k)/\sO_{\sM}$ such
that any local lifting is
a heat operator.
The {\it symbol map} of $\rm H$ is
$\epsilon\circ {\rm H}:f^*T_{\sM}\to S^2T_{S/\sM}$.
A heat operator ${\rm H}$
induces an $\sO_{\sM}$-map,
denoted by the same $\rm H$,
$T_\sM\to f_*\sW(\sF)$
The heat operator and the preceding induced map
will be used interchangeably throughout.
A (projective) heat operator on $\Theta^k$
determines a (projective) connection
$\nabla^{\rm H}$
on $f_*\Theta^k$ in a natural way \cite{GJ}.
\end{defn}

We recall firstly a general result (forget moduli) of G. Faltings 
(cf. \cite{F}).
Let $Z\to S$ be smooth and $\sK=K_{Z/S}$. Consider
$$0\to\sD^{\le 1}_{Z/S}(\sL)\to 
\sD^{\le 2}_{Z/S}(\sL)\xrightarrow{\e_2}
S^2T_{Z/S}\to 0$$
$$0\to\sD^{\le 1}_{Z/S}(\sL)\to 
S^2\sD^{\le 1}_{Z/S}(\sL)\xrightarrow{S^2\e_1}
S^2T_{Z/S}\to 0.$$
For any $\rho\in H^0(Z,S^2T_{Z/S})$, let $a(\sL,\rho)$, $b(\sL,\rho)$
be the obstruction classes to lift $\rho$ to 
$H^0(Z,\sD^{\le 2}_{Z/S}(\sL))$ and to $H^0(Z,S^2\sD^{\le 1}_{Z/S}(\sL))$
respectively. Then 

\begin{prop}\label{prop0.0} Under 
$\sD^{\le1}_{Z/S}(\sL^k)\cong\sD^{\le1}_{Z/S}(\sL)$, one has
\item i) $b(\sL^k,\rho)=kb(\sL,\rho)$ for any $k\in\Bbb Q$.
\item ii) $a(\sO_Z,\rho)\in H^1(\sO_Z)\oplus
H^1(T_{Z/S})$ has zero projection in $H^1(\sO_Z)$.
\item iii) There is a class $c(\sK,\rho)\in H^1(\sO_Z)$, independent of $\sL$,
such that $$2a(\sL,\rho)=b(\sL,\rho)+\,
^tb(\sL^{-1}\otimes\sK,\rho)+c(\sK,\rho).$$
\end{prop}

We come back to construct the projective connection. Recall that
$$\CD
X=\sC\times_{\sM}S@>\pi>> S\\
@VVV         @V f VV\\
\sC@>>> \sM
\endCD$$
For simplicity, we assume that
there exists a universal bundle $E$ on $X$, let $\E=\sE nd^0(E)$. 
We will arrive at a subsheaf $\sG_E\subset\,^{tr}\sA_E^{-1}$ 
(see \cite{BS} for details of $\,^{tr}\sA_E^{-1}$) fitting into an exact
sequence 
\begin{gather}\label{eq1.1}
0\to\omega_{X/S}\to\sG_E\xrightarrow{\rm res}\E\to 0, 
\end{gather}
which induces, by taking 2-th symmetric tensor, the exact sequence 
$$0\to\sG_E\to S^2(\sG_E)\otimes T_{X/S}
\xrightarrow{\Sym^2(\rm res)\otimes \id}
S^2(\E)\otimes T_{X/S}\to 0.$$ 
Define $S(\sG_E):=(\Sym^2({\rm res})\otimes \id)^{-1}(\id\otimes T_{X/S})$, 
which fits into 
\begin{gather}\label{eq1.2}
0\to\sG_E\xrightarrow{\iota} S(\sG_E)\xrightarrow{q}T_{X/S}\to 0,
\end{gather} 
where $q=\Sym^2({\rm res})\otimes \id$, 
$\iota(\alpha)=\Sym^2(\alpha\otimes dt)\otimes\partial_t$ locally. Let
\begin{gather}\label{eq0.1}
0\to R^1\pi_*(\omega_{X/S})\to R^1\pi_*(\sG_E)\xrightarrow{res} R^1\pi_*(\E)\to0
\end{gather}
\begin{gather}\label{eq0.2}
0\to R^1\pi_*(\sG_E)\xrightarrow{\iota}R^1\pi_*S(\sG_E)\xrightarrow{q} R^1\pi_*(T_{X/S})\to 0 
\end{gather}
be the exact sequences induced by (\ref{eq1.1}), (\ref{eq1.2}).
Let $\Delta\subset X\times_SX$ be the diagonal, consider the induced
diagram
$$\CD
0 @>>>\sG_E\boxtimes\sG_E @>>>\sG_E\boxtimes\sG_E(\Delta)@>>>
\sG_E\boxtimes\sG_E(\Delta)|_{\Delta} @>>>0   \\
@.
@VVV             @VVV          @VVV        @.\\
0@>>>\E\boxtimes \E@>>>\E\boxtimes \E(\Delta)@>>>
\E\boxtimes \E(\Delta)|_{\Delta}@>>> 0 
\endCD$$
on $X\times_SX$ ($\E\boxtimes \E$ denotes $p_1^*\E\otimes p_2^*\E$).
All vertical maps are induced by $\sG_E\xrightarrow{res}\E\to 0$, thus 
(\ref{eq1.2}) is a sub-sequence of the
rightmost vertical map above.  Let $\sF_1\subset\sG_E\boxtimes\sG_E(\Delta)$,
$\sF_2\subset \E\boxtimes \E(\Delta)$ be subsheaves
satisfying 
$$\CD
0 @>>>\sG_E\boxtimes\sG_E @>>>\sF_1@>>>
S(\sG_E) @>>>0   \\
@. @VVV             @VVV          @V{q}VV        @.\\
0@>>>\E\boxtimes \E@>>>\sF_2@>>>
T_{X/S}@>>> 0 
\endCD$$

Taking direct images, considering the connecting maps, we have
$$\CD
R^1\pi_*S(\sG_E) @>\tilde\delta>> S^2R^1\pi_*(\sG_E)\\
@V q VV                        @V S^2(res)VV \\  
R^1\pi_*(T_{X/S}) @>\delta>> S^2R^1\pi_*(\E)
\endCD$$
which induces the commutative diagram
$$\CD
R^1\pi_*(\sG_E)@>\delta_1>>R^1\pi_*(\sG_E)\\
  @V\iota VV             @V\Sym^2VV  \\
R^1\pi_*S(\sG_E) @>\tilde\delta>>S^2R^1\pi_*(\sG_E)\\
@V q VV               @V S^2(res)VV        \\
R^1\pi_*(T_{X/S}) @>\delta>> S^2R^1\pi_*(\E)
\endCD$$
where $\delta_1$ is defined such that the diagram is commutative,
the first vertical exact sequence is
(\ref{eq0.2}), 
the second vertical exact sequence is induced
by taking 2-th symmetric tensor of (\ref{eq0.1}) (note that
$\sO_S\cong R^1\pi_*(\omega_{X/S})$).  
For $\Sym^2$ see Remark \ref{prop0.2r}. 

\begin{lem}\label{prop0} The map $\delta_1:R^1\pi_*(\sG_E)\to R^1\pi_*(\sG_E)$ 
is the identity map.\end{lem}

\begin{proof} See Lemma \ref{prop5.13}.\end{proof}

\begin{thm}\label{prop0.2} i) There is an isomorhism 
$\phi:R^1\pi_*(\sG_E)\cong \sD^{\le 1}_{S/\sM}(\lambda_{\E})$ such that
the following diagram is commutative
$$\CD
0@>>>\sO_S=R^1\pi_*\omega_{X/S}@>>>
R^1\pi_*(\sG_E)@>{\rm res}>>
R^1\pi_*(\E)@>>>0\\
@. @V(2r)\cdot{\rm id}VV     @V\phi VV @V\vartheta VV  @.  \\
0@>>>\sO_S=R^0\pi_*\omega_{X/S}[1]@>>>
\sD^{\le 1}_{S/\sM}(\lambda_{\E})
@>\e_1>>T_{S/\sM}@>>>0
\endCD$$

ii) For any affine covering $\{\sM_i\}_{i\in I}$ of $\sM$, on each
$S_i:=f^{-1}(\sM_i)$, there is an isomorphism
$\phi_i:R^1\pi_*S(\sG_E)\cong \sD^{\le 1}_S(\lambda_{\E})$
such that 
$$\CD
0@>>>R^1\pi_*(\sG_E)@>2\cdot\iota>>R^1\pi_*S(\sG_E)@>q>>R^1\pi_*(T_{X/S})@>>>0\\
@.@V\phi VV              @V\phi_i VV @VVV @.\\
0@>>>\sD^{\le 1}_{S/\sM}(\lambda_{\E})@>>>\sD^{\le 1}_S(\lambda_{\E})
@>\sigma>>f^*T_{\sM}@>>>0
\endCD$$
is commutative on $S_i$. Moreover, $\{\phi_i-\phi_j\}$ define a class
in ${\rm H}^1(\Omega^1_{\sM})$.
\end{thm}

\begin{proof} See Theorem \ref{prop5.4}, Corollary \ref{prop5.9}
and Remark \ref{rem3.18}.
\end{proof}

By the above theorem, on each $S_i$, we get commutative diagram
$$\CD
\sD^{\le 1}_{S/\sM}(\lambda_{\E})@>\phi^{-1}>>R^1\pi_*(\sG_E)@>\delta_1>>
R^1\pi_*(\sG_E)@>2r\cdot\phi>>\sD^{\le 1}_{S/\sM}(\lambda_{\E})\\
@VVV  @V2\cdot\iota VV             @V\Sym^2VV @VVV \\
\sD^{\le 1}_S(\lambda_{\E})@>\phi_i^{-1}>>R^1\pi_*S(\sG_E) @>\tilde\delta>> 
S^2R^1\pi_*(\sG_E)@>S^2(\phi)>>S^2\sD^{\le 1}_{S/\sM}(\lambda_{\E})\\
@V\sigma VV @VVV               @V q VV        @VS^2\e_1VV\\
f^*T_{\sM}@>\cong>> R^1\pi_*(T_{X/S}) @>\delta>> S^2R^1\pi_*(\E)@>\cong>> 
S^2T_{S/\sM}.
\endCD$$
\begin{remark}\label{prop0.2r}
i) For a precise definition of ${\rm Sym}^2$ above   
we refer to proof of Lemma \ref{prop5.13}; ii)  
The insertion of $2r$ in $2r\cdot\phi$ in the 1st row 
is due to Theorem \ref{prop0.2} i)
combined with the definition of ${\rm Sym}^2$.      
\end{remark}
The above diagram gives the following commutaive diagram on $S_i$    
$$\CD
0@>>>\sD^{\le 1}_{S/\sM}(\lambda_{\E})@>>>\sD^{\le 1}_S(\lambda_{\E})@>\sigma>>
f^*T_{\sM}@>>> 0 \\
@.  @|    @V\tilde\delta^i_{\rm H} VV  @V\delta_{\rm H} VV   @. \\   
0@>>>\sD^{\le 1}_{S/\sM}(\lambda_{\E})@>>>S^2\sD^{\le 1}_{S/\sM}(\lambda_{\E})
@>S^2\e_1>> S^2T_{S/\sM}@>>> 0 
\endCD$$
where $\delta_{\rm H}$, $\tilde\delta^i_{\rm H}$ are defined in a clear way
such that each $\tilde\delta^i_{\rm H}$ 
induces an identity map on $\sD^{\le 1}_{S/\sM}(\lambda_{\E})$.

\begin{prop}\label{prop0.3} For any 
$\rho=\delta_{\rm H}(v)\in H^0(S, S^2T_{S/\sM})$,
where $v\in H^0(S,f^*T_{\sM})$ and $\sM$ is replaces by its affine open
set, one has
\item i) $2a(\sL,\rho)=b(\sL,\rho)+\,^tb(\sL^{-1}\otimes K_{S/\sM},\rho)$.
\item ii) When $\sL=K_{S/\sM}^{\mu}$, where $\mu\in\Bbb Q$ and $\mu\neq 1$,
one has $$a(K_{S/\sM}^{\mu},\rho)=
\frac{2\mu-1}{2\mu}b(K^{\mu}_{S/\sM},\rho).$$ 
\end{prop}

\begin{proof} To prove i), by Proposition \ref{prop0.0} iii), it is enough to show that
$c(\sK,\rho)=0$. The class is independent of $\sL$. Thus, by 
taking $\sL=\sO_S$ and using Proposition \ref{prop0.0} ii), it is enough 
to show that 
$$\,^tb(\sK,\rho)\in H^1(\sD^{\le1}_{S/\sM}(\sO_S))=
H^1(\sO_S)\oplus H^1(T_{S/\sM})$$
has trivial projection in $H^1(\sO_S)$ where  
$\sK=K_{S/\sM}$.  
We have (noting $\sK=\lambda_{\E}$) 
$$\CD
\sD^{\le 1}_{S/\sM}(\sK)@=\sD^{\le 1}_{S/\sM}(\sK)\\
@VVV               @VVV \\
\sD^{\le 1}_S(\sK)@>\tilde\delta_{\rm H}>>S^2\sD^{\le 1}_{S/\sM}(\sK)\\
@V\sigma VV                       @V S^2\e_1VV\\
f^*T_{\sM} @>\delta_{\rm H}>>S^2T_{S/\sM}
\endCD$$
Let $\{\sU_i\}_{i\in I}$ be an affine covering of $S$ and 
$v_i\in T_S(\sU_i)$ be local liftings of $v\in f^*T_{\sM}(S)$. Let
$\{d_i\in\sD^{\le1}_S(\sK)(\sU_i)\}_{i\in I}$ be 
such that $\e_1(d_i)=v_i$ ($i\in I$). Then, by above diagram, 
$b(\sK,\rho)=
\{d_i-d_j\in\sD^{\le1}_{S/\sM}(\sK)(\sU_i\cap\sU_j)\}.$
Thus the class $\,^tb(\sK,\rho)\in H^1(\sD^{\le1}_{S/\sM}(\sO_S))=
H^1(\sO_S)\oplus H^1(T_{S/\sM})$ is defined by
the cocycle $\{\,^td_i-\,^td_j\}$, where 
$\,^td_i:=A_i-v_i\in \sD^{\le1}_S(\sO_S)(\sU_i).$
Thus the projection of $\,^tb(\sK,\rho)$ in $H^1(\sO_S)$ is defined
by $\{A_i-A_j\}$, which is a trivial class.

To show ii), we remark that for any nonzero $\mu\in\Bbb Q$, through 
canonical isomorphisms 
$\psi_{\mu}: \sD^{\le1}(\sK)\cong\sD^{\le1}(\sK^{\mu})$,
the above diagram induces
$$\CD
\sD^{\le 1}_{S/\sM}(\sK^{\mu})@>\bullet\mu>>\sD^{\le 1}_{S/\sM}(\sK^{\mu})\\
@VVV               @VVV \\
\sD^{\le 1}_S(\sK^{\mu})@>\tilde\delta^{\mu}_{\rm H}>>
S^2\sD^{\le 1}_{S/\sM}(\sK^{\mu})\\
@V\sigma VV                       @V S^2\e_1VV\\
f^*T_{\sM} @>\delta_{\rm H}>>S^2T_{S/\sM}
\endCD$$
where $\tilde\delta^{\mu}_{\rm H}=
S^2\psi_{\mu}\circ\tilde\delta_{\rm H}\circ\psi^{-1}_{\mu}$. 
Using the above diagram, 
we can compute $\,^tb(\sK^{1-\mu},\rho)$. Let
$\{d_i\in\sD^{\le1}_S(\sK^{1-\mu})(\sU_i)\}_{i\in I}$ be 
such that $\e_1(d_i)=v_i$ ($i\in I$). Then
$b(\sK^{1-\mu},\rho)=(1-\mu)\{d_i-d_j\}$,
which implies that 
$$\,^tb(\sK^{1-\mu},\rho)=(1-\mu)\{\,^td_i-\,^td_j\}.$$
On the other hand, 
$\{-\,^td_i\in\sD^{\le1}_S(\sK^{\mu})(\sU_i)\}_{i\in I}$ are local
liftings of $v$, which means that
$b(\sK^{\mu},\rho)=-\mu\{\,^td_i-\,^td_j\}=
\frac{\mu}{\mu-1}\,^tb(\sK^{1-\mu},\rho)$. Thus 
$$a(\sK^{\mu},\rho)=\frac{2\mu-1}{2\mu}b(\sK^{\mu},\rho).$$
\end{proof}

\begin{thm}\label{prop0.4} Replace $\sM$ by its affine open sets, let $\{\sU_i\}_{i\in I}$ be an
affine open covering of $S$. Then, for any $v\in f^*T_{\sM}(S)$, there
are   
$$d^i_S\in \sD^{\le1}_S(\sK^{\mu})(\sU_i),\quad 
d^i_{S/\sM}\in \sD^{\le1}_{S/\sM}(\sK^{\mu})(\sU_i),\quad
D^i_{S/\sM}\in \sD^{\le2}_{S/\sM}(\sK^{\mu})(\sU_i),$$
where $\sK:=K_{S/\sM}$, such that 
$$\left\{{\rm H}(v)_i:=d^i_S-d^i_{S/\sM}+\frac{2}{1-2\mu}D^i_{S/\sM}
\in\sD^{\le2}_S(\sK^{\mu})(\sU_i)\right\}_{i\in I}$$
form a global section ${\rm H}(v)\in H^0(S,\sD^{\le2}_S(\sK^{\mu}))$ with
$$\sigma({\rm H}(v))=v,\quad \e_2({\rm H}(v))=
\frac{2}{1-2\mu}\delta_{\rm H}(v).$$
\end{thm}

\begin{proof} Let
$\{d^i_S\in\sD^{\le 1}_S(\sK^{\mu})(\sU_i)\}_{i\in I}$ be such that 
$\sigma(d^i_S)=v|_{\sU_i}$.
Then $$\{\mu(d^i_S-d^j_S)\in\sD^{\le 1}_{S/\sM}(\sK^{\mu})(\sU_i\cap \sU_j)\}$$
defines the class $b(\sK^{\mu},\delta_{\rm H}(v))\in 
H^1(S,\sD^{\le 1}_{S/\sM}(\sK^{\mu}))$, which is the obstruction
for lifting $\delta_{\rm H}(v)\in H^0(S,S^2T_{S/\sM})$ to 
$H^0(S,S^2\sD^{\le 1}_{S/\sM}(\sK^{\mu}))$. Let 
$$\{D^i_{S/\sM}\in \sD^{\le2}_{S/\sM}(\sK^{\mu})(\sU_i)\}_{i\in I}$$
be local liftings of $\delta_{\rm H}(v)$. Then, 
by Proposition \ref{prop0.3}, 
$$\{d^i_S-d^j_S\}=\frac{2}{2\mu-1}\{D^i_{S/\sM}-D^j_{S/\sM}\}$$
as cohomology classes. Thus there are 
$\{d^i_{S/\sM}\in \sD^{\le1}_{S/\sM}(\sK^{\mu})(\sU_i)\}_{i\in I}$
satisfying the requirements in the theorem.
\end{proof} 

\begin{cor}\label{prop0.5}
There exists uniquely a projective heat operator, 
$${\rm H}: T_{\sM}\to f_*\sW(\Theta^k)/\sO_{\sM}$$
such that $(f_*\e_2)\cdot{\rm H}:T_{\sM}\to f_*S^2T_{S/\sM}$ coincides with
$f_*\delta_{\rm H}$.
\end{cor}

\begin{proof} For any open set $U\subset\sM$ and $v\in T_{\sM}(U)$,
by Theorem \ref{prop0.4}, we can construct 
a ${\rm H}(v)\in f_*\sW(\Theta^k)(U)$.
If ${\rm H}(v)'$ is another such operator, ${\rm H}(v)-{\rm H}(v)'$
must have symbol in $H^0(f^{-1}(U),T_{S/\sM})=0$, so 
$${\rm H}(v)-{\rm H}(v)'\in H^0(f^{-1}(U),\sO_S)=f_*\sO_S(U)=\sO_{\sM}(U).$$
Hence a unique map $T_{\sM}\to f_*\sW(\Theta^k)/\sO_{\sM}$.
\end{proof}

Now we construct the logarithmic extension of above operator. Let  
$\ov\sC\to \ov\sM$ be the family in Introduction. 
Let $\sU(r,d)\to\ov\sM$ be the family
of moduli spaces of semistable (for canonical polarization) torsion free
sheaves of rank $r$ and degree $d$. Fix a line bundle $\sN$ on $\ov\sC$
of relative degree $d$, let $f:S\to\sM$ be the family
of moduli spaces of stable bundles of rank $r$ with fixed
determinant $\sN|_C$ ($[C]\in\sM$). 

\begin{notation}\label{prop2}
Let 
$f_Z:Z\to\ov\sM$
be defined as the Zariski closure of $S$ in $\sU(r,d)$ and
$f_T:T\to\ov\sM$
be the open set of $Z$ consisting of locally free sheaves. Let 
$\ov f:\ov S\subset T\to\ov\sM$ 
be the open set of stable bundles.
Then $\ov f: \ov S\to\ov\sM$ is smooth (cf. Lemma \ref{prop9.9}). 
\end{notation}

Let $E$ be the unversal bundle on $\ov X$, where $\ov X$ is defined by
$$\CD
\ov X=\ov\sC\times_{\ov\sM}\ov S@>\pi>>\ov S\\
@VVV         @V \ov f VV\\
\ov\sC@>>> \ov\sM
\endCD$$
Let $D\subset\ov X$ 
be the divisor of singular curves. Then we have 
$$\sG_E\subset\, ^{tr}\sA^{-1}_E({\rm log}D)$$ 
fitting into the exact sequence
\begin{gather}\label{eq1.3}
0\to\omega_{\ov X/\ov S}\to\sG_E\xrightarrow{res}\E:=\sE nd^0(E)\to 0. 
\end{gather}
Similarly, there is a sheaf $S(\sG_E)\subset {\rm Sym}^2(\sG_E)\otimes
\omega^{-1}_{\ov X/\ov S}$ fitting into
\begin{gather}\label{eq0.3}
0\to\sG_E\xrightarrow{\iota} S(\sG_E)\xrightarrow{q}\omega^{-1}_{\ov X/\ov S}\to 0.
\end{gather}
They induce the following exact sequences
\begin{gather}\label{eq0.4}
0\to R^1\pi_*(\omega_{\ov X/\ov S})\to R^1\pi_*(\sG_E)\xrightarrow{res}
R^1\pi_*(\E)\to 0,
\end{gather}
\begin{gather}\label{eq0.5}
0 \to R^1\pi_*(\sG_E)\xrightarrow{\iota} R^1\pi_*S(\sG_E)\xrightarrow{q} 
R^1\pi_*(\omega^{-1}_{\ov X/\ov S})\to 0.
\end{gather}

Let $\ov\Delta\subset \ov X\times_{\ov S}\ov X:=P$ be the diagonal and
$\sO(\ov\Delta)$ be the dual of its ideal sheaf. Then
$0\to\omega_{P/\ov S}\to\omega_{P/\ov S}(\ov\Delta)\to
\sE xt_P^1(\sO_{\ov\Delta},\omega_{P/\ov S})\to 0,$
one checks that the relative dualizing sheaf $\omega_{P/\ov S}$ is
$\omega_{\ov X/\ov S}\boxtimes\omega_{\ov X/\ov S}$ and 
$\sE xt_P^1(\sO_{\ov\Delta},\omega_{P/\ov S})$ is the relative dualizing
sheaf of $\ov\Delta/\ov S$. Thus we have
\begin{gather}\label{eq0.6}
0\to\sO\to\sO(\ov\Delta)\to\omega^{-1}_{\ov X/\ov S}\to 0,
\end{gather}
which similarly induces the commutative diagram
$$\CD
R^1\pi_*(\sG_E)@>\delta_1>>R^1\pi_*(\sG_E)\\
  @V\iota VV             @V\Sym^2VV  \\
R^1\pi_*S(\sG_E) @>\tilde\delta>>S^2R^1\pi_*(\sG_E)\\
@V q VV               @V S^2(res)VV        \\
R^1\pi_*(\omega^{-1}_{\ov X/\ov S}) @>\delta>> S^2R^1\pi_*(\E)
\endCD$$
where $\delta_1$ denote the map induced by $\tilde\delta$,
the first vertical exact sequence is
(\ref{eq0.5}), 
the second vertical exact sequence is induced
by taking 2-th symmetric tensor of (\ref{eq0.4}) (note that
$\sO_{\ov S}\cong R^1\pi_*(\omega_{\ov X/\ov S})$). 

Let $B=\ov\sM\setminus\sM$ and $W={\ov f}^{-1}(B)\subset \ov S$. Consider
$$0\to T_{\ov S/\ov\sM}\to T_{\ov S}\xrightarrow{d{\ov f}}
{\ov f}^*T_{\ov\sM}\to 0,$$
$$0\to\sD^{\le1}_{\ov S/\ov\sM}(\sL)\to\sD^{\le1}_{\ov S}(\sL)
\xrightarrow{\sigma}{\ov f}^*T_{\ov\sM}\to 0.$$

\begin{notation} Let $T_{\ov\sM}(B)\subset T_{\sM}$ be the subsheaf of
vector fields that preserve $B$. Let $T_{\ov S}(\log W)\subset T_{\ov S}$,
$\sD^{\le 1}_{\ov S}(\sL)(\log W)
\subset\sD^{\le 1}_{\ov S}(\sL)$ 
be the subsheaves such that the following are exact sequences
\begin{gather}\label{0.6}
0\to T_{\ov S/\ov\sM}\to T_{\ov S}(\log W)\xrightarrow{d{\ov f}}
{\ov f}^*T_{\ov\sM}(B)\to 0,
\end{gather}
\begin{gather}\label{07}
0\to\sD^{\le1}_{\ov S/\ov\sM}(\sL)\to
\sD^{\le1}_{\ov S}(\sL)(\log W)
\xrightarrow{\sigma}{\ov f}^*T_{\ov\sM}(B)\to 0.
\end{gather}
\end{notation}

\begin{lem}\label{prop1} The map $\delta_1:R^1\pi_*(\sG_E)\to R^1\pi_*(\sG_E)$ is
identity.\end{lem}

\begin{proof} Since $S\subset\ov S$ is dense in $\ov S$ by definition
(cf. Notation \ref{prop2}), the lemma follows from Lemma \ref{prop0}.
\end{proof}

\begin{thm}\label{prop0.6} i) There is an isomorhism 
$\phi:R^1\pi_*(\sG_E)\cong \sD^{\le 1}_{\ov S/\ov \sM}(\lambda_{\E})$ 
such that
the following diagram is commutative
$$\CD
0@>>>\sO_{\ov S}=R^1\pi_*\omega_{\ov X/\ov S}@>>>
R^1\pi_*(\sG_E)@>{\rm res}>>
R^1\pi_*(\E)@>>>0\\
@. @V(2r)\cdot{\rm id}VV     @V\phi VV @V\vartheta VV  @.  \\
0@>>>\sO_{\ov S}=R^0\pi_*\omega_{\ov X/\ov S}[1]@>>>
\sD^{\le 1}_{\ov S/\ov\sM}(\lambda_{\E})
@>\e_1>>T_{\ov S/\ov\sM}@>>>0
\endCD$$

ii) For any affine covering $\{\ov\sM_i\}_{i\in I}$ of $\ov\sM$, on each
$\ov S_i:={\ov f}^{-1}(\ov\sM_i)$, there is an isomorphism
$\phi_i:R^1\pi_*S(\sG_E)\cong \sD^{\le 1}_{\ov S}(\lambda_{\E})(\log W)$
such that 
$$\CD
0@>>>R^1\pi_*(\sG_E)@>2\cdot\iota>>R^1\pi_*S(\sG_E)@>q>>
R^1\pi_*(\omega^{-1}_{\ov X/\ov S})@>>>0\\
@.@V\phi VV              @V\phi_i VV @VVV @.\\
0@>>>\sD^{\le 1}_{\ov S/\ov\sM}(\lambda_{\E})@>>>
\sD^{\le 1}_{\ov S}(\lambda_{\E})(\log W)
@>\sigma>>{\ov f}^*T_{\ov\sM}(B)@>>>0
\endCD$$
is commutative on $\ov S_i$, and $\{\phi_i-\phi_j\}$ define a class
in ${\rm H}^1(\Omega^1_{\ov\sM}(\log B)$.
\end{thm}

\begin{proof} See Proposition \ref{prop6.4} and
Theorem \ref{prop6.6}.
\end{proof}

Similarly, we have the commutative diagram on each $\ov S_i$    
$$\CD
0@>>>\sD^{\le 1}_{\ov S/\ov\sM}(\lambda_{\E})@>>>
\sD^{\le 1}_{\ov S}(\lambda_{\E})(\log W)@>\sigma>>
{\ov f}^*T_{\ov\sM}(B)@>>> 0 \\
@.  @|    @V\tilde\delta^i_{\rm H} VV  @V\delta_{\rm H} VV   @. \\   
0@>>>\sD^{\le 1}_{\ov S/\ov\sM}(\lambda_{\E})@>>>
S^2\sD^{\le 1}_{\ov S/\ov\sM}(\lambda_{\E})
@>S^2\e_1>> S^2T_{\ov S/\ov\sM}@>>> 0 
\endCD$$

\begin{prop}\label{prop0.7} For any $\rho=\delta_{\rm H}(v)\in 
H^0(\ov S, S^2T_{\ov S/\ov\sM})$, where $\ov\sM$ is replaced by its
affine open sets and
$v\in H^0(\ov S,{\ov f}^*T_{\ov\sM}(B))$, one has
\item i)  $2a(\sL,\rho)=b(\sL,\rho)+\,^tb(\sL^{-1}\otimes K_{\ov S/\ov\sM},\rho)$.
\item ii) When $\sL=K_{\ov S/\ov\sM}^{\mu}$, where $\mu\in\Bbb Q$ and $\mu\neq 1$,
one has $$a(K_{\ov S/\ov\sM}^{\mu},\rho)=
\frac{2\mu-1}{2\mu}b(K^{\mu}_{\ov S/\ov\sM},\rho).$$ 
\end{prop}

\begin{proof}  Note that we still have 
$\sK=K_{\ov S/\ov\sM}=\lambda_{\E}$, the proof is the same as 
that of Proposition \ref{prop0.3}. We just remark two points:
(1) for any operator $d$ of 
$\sD^{\le 1}_{\ov S}(\sL)(\log W)
\subset\sD^{\le 1}_{\ov S}(\sL)$, its adjoint operator $\,^td$ is still
in $\sD^{\le 1}_{\ov S}(\sL^{-1}\otimes\sK)(\log W)$.
(2) for any nonzero $\mu\in\Bbb Q$, the 
canonical isomorphism
$\psi_{\mu}: \sD^{\le 1}_{\ov S}(\sK)
\cong\sD^{\le}_{\ov S}(\sK^{\mu})$
induces an isomorphism 
$\sD_{\ov S}^{\le1}(\sK)(\log W)\cong\sD_{\ov S}^{\le1}(\sK^{\mu})(\log W)$.
\end{proof}

\begin{thm}\label{prop0.8} Replace $\ov\sM$ by its affine open set, 
let $\{\sU_i\}_{i\in I}$ be an
affine open covering of $\ov S$. Then, for any 
$v\in {\ov f}^*T_{\ov\sM}(B)(\ov S)$, there
are 
$$d^i_{\ov S}\in \sD^{\le1}_{\ov S}(\sK^{\mu})(\log W)(\sU_i),\quad 
d^i_{\ov S/\ov\sM}\in \sD^{\le1}_{\ov S/\ov\sM}(\sK^{\mu})(\sU_i),$$
and $D^i_{\ov S/\ov\sM}\in \sD^{\le2}_{\ov S/\ov\sM}(\sK^{\mu})(\sU_i)$
such that $$\left\{{\rm H}(v)_i:=d^i_{\ov S}-d^i_{\ov S/\ov\sM}+
\frac{2}{1-2\mu}D^i_{\ov S/\ov\sM}
\in\sD^{\le2}_{\ov S}(\sK^{\mu})(\sU_i)\right\}_{i\in I}$$
form a global section ${\rm H}(v)\in H^0(\ov S,\sD^{\le2}_{\ov S}(\sK^{\mu}))$ with
$$\sigma({\rm H}(v))=v,\quad \e_2({\rm H}(v))=
\frac{2}{1-2\mu}\delta_{\rm H}(v).$$
\end{thm}

\begin{proof} Let
$\{d^i_{\ov S}\in\sD^{\le 1}_{\ov S}(\sK^{\mu})(\log W)(\sU_i)\}_{i\in I}$ 
be such that 
$\sigma(d^i_{\ov S})=v|_{\sU_i}$.
Then $$\{\mu(d^i_{\ov S}-d^j_{\ov S})\in\sD^{\le 1}_{\ov S/\ov\sM}(\sK^{\mu})(\sU_i\cap \sU_j)\}$$
defines the class $b(\sK^{\mu},\delta_{\rm H}(v))\in 
H^1(\ov S,\sD^{\le 1}_{\ov S/\ov\sM}(\sK^{\mu}))$, which is the obstruction
for lifting $\delta_{\rm H}(v)\in H^0(\ov S,S^2T_{\ov S/\ov\sM})$ to 
$H^0(\ov S,S^2\sD^{\le 1}_{\ov S/\ov\sM}(\sK^{\mu})$. Let 
$$\{D^i_{\ov S/\ov\sM}\in \sD^{\le2}_{\ov S/\ov\sM}(\sK^{\mu})(\sU_i)\}_{i\in I}$$
be local liftings of $\delta_{\rm H}(v)$. Then, 
by Proposition \ref{prop0.2}, 
$$\{d^i_{\ov S}-d^j_{\ov S}\}=
\frac{2}{2\mu-1}\{D^i_{\ov S/\ov\sM}-D^j_{\ov S/\ov\sM}\}$$
as cohomology classes. Thus there are 
$\{d^i_{\ov S/\ov\sM}\in \sD^{\le1}_{\ov S/\ov\sM}(\sK^{\mu})(\sU_i)\}_{i\in I}$
satisfying the requirements in the theorem.
\end{proof} 

\begin{cor}\label{prop0.9} There exists uniquely a projective heat operator
$$\ov{\rm H}: T_{\ov\sM}(B)\to 
\ov f_*\sD^{\le 2}_{\ov S}(\Theta^k)_B/\sO_{\ov\sM}$$
such that $(\ov f_*\e_2)\cdot\ov{\rm H}:T_{\ov\sM}(B)\to f_*S^2T_{S/\sM}$ 
coincides with $\ov f_*\delta_{\rm H}$. Moreover, $\ov{f}_*(\Theta^k)$
is a coherent sheaf on $\ov\sM$. 
\end{cor}

\begin{proof} The coherence of  
$\ov{f}_*(\Theta^k)$ follows from Theorem \ref{prop9.2}.
Since $\ov {f}_*T_{\ov S/\ov\sM}=0$ is still true, the rest follows
the same proof of
Corollary \ref{prop0.5} if 
$\ov f_*\sO_{\ov S}=\sO_{\ov\sM}$. 
We will prove (cf. Proposition \ref{prop9.1}) that 
the fibres of $\ov f:\ov S\to\ov\sM$ is
dense in the fibres of $f_Z:Z\to\ov\sM$. Thus $Z\setminus\ov S$ 
has codimension
at least $2$. By passing to the normalization $\iota:\widetilde{Z}\to Z$,
$\iota^{-1}(\ov S)\cong\ov S$ and 
${\rm codim}(\widetilde Z\setminus \iota^{-1}(\ov S))\ge 2$ since $\ov S$ is the
open set of smooth points of $Z$ and $\iota$ is finite. We have 
$\ov f_*\sO_{\ov S}=(\ov f\cdot\iota)_*\sO_{\widetilde Z}=\sO_{\ov\sM}.$
\end{proof}

\section{First Order Differential Operators of the Determinant
Bundle}

We give at first a short review of what
we need from \cite{BS}. 
Let $\pi:X\to S$ be a smooth proper morphism of relative
dimension $1$ between smooth varieties in
characteristic 0.  We write $K_{X/S}$ or $\omega_{X/S}$
interchangeably for the dualizing sheaf.
One has an exact sequence
\begin{gather}\label{eq4.1}
0\to T_{X/S}\to T_X \xrightarrow{d\pi} \pi^*T_S\to 0.
\end{gather}
As in \cite{BS}, one defines the subsheaf $\pi^{-1}T_S\subset \pi^*T_S$
and its preimage $T_{\pi}=d\pi^{-1}T_\pi\subset T_X$, defining
the  exact sequence
\begin{gather}\label{eq4.2}
0\to T_{X/S}\to T_{\pi}\xrightarrow{d\pi} \pi^{-1}T_S\to 0.
\end{gather}
Let $E$ be a vector bundle on $X$, and
$\lambda_E=\det R\pi_*E$ be its determinant
bundle.
The Atiyah algebra $\sA_E$ is the subalgebra
of the sheaf of first order differential operators
on $E$ with symbolic part in $({\rm id}\otimes T_X)\cong T_X$.
The relative Atiyah algebra $\sA_{E/S} \subset\sA_{E}$
consists of those differential operators with
symbolic part in $T_{X/S}$, and
$\sA_{E,\pi}\subset \sA_E$ with symbolic part in $T_{\pi}$.
Let  $^{tr}\sA_E^{-1}$ be the subquotient of the sheaf
defined in \cite{BS}
$$E\boxtimes_{\sO_S} (E^*\otimes \omega_{X/S})(2 \Delta)/E\boxtimes_{\sO_S}
(E^*\otimes \omega_{X/S})(-\Delta)$$
where $\Delta\subset X\times_S X$ denotes the diagonal,
which fits into an exact sequence
\begin{gather}\label{eq4.3}
0\to \omega_{X/S}\to \,^{tr}\sA_E^{-1}\buildrel{\rm res}\over
\longrightarrow \sA_{E/S}\to 0.
\end{gather}
The trace complex is defined by
\begin{gather} \label{eq4.4}
^{tr}\sA_E^{\bullet}:
\sO_X \xrightarrow{d_{X/S}}  {^{tr}\sA_E^{-1}}
\buildrel{\rm res}\over \longrightarrow \sA_{E,\pi}
\end{gather}
with $\sA_{E,\pi} $ in degree 0.
One has
\begin{prop}\label{prop4.1}
$^{tr}\sA_E^{\bullet}$ carries an algebra structure for which
$R^0\pi_*(^{tr}\sA_E^{\bullet})$ is canonically isomorphic to
$\sA_{\lambda_E}$ (\cite{BS}, 2.3.1, see also \cite{ET}).
\end{prop}
For the purpose of this paper it is more convient to
define the trace complex concentrated
only on $i=-1$ and $i=0$ of the original
trace complex.  This modified trace complex
is still denoted by $^{tr}\sA_E^{\bullet}$ whose
$0$-th direct image is easily seen to be the same
as that of the orginal one.
With the modified trace complex,
one has now an exact sequence
\begin{gather}\label{eq4.7}
0\to \omega_{X/S}[1]\to \,^{tr}\sA_{E}^{\bullet}\to \sA^{\bullet}_{E,\pi}\to 
0,
\end{gather}
where the complex $\sA^{\bullet}_{E,\pi}$ is defined by
\begin{gather} \label{eq4.8}
\sA^{\bullet}_{E,\pi}: \sA_{E/S} \to \sA_{E,\pi},
\end{gather}
and thus is quasi-isomorphic to $\pi^{-1}T_S$.

\begin{notation} \label{not:rel}
Let $\pi:X\to S$ be as before, and
$f:S\to M$ be a smooth morphism where
$M$ is a smooth variety.
Denote by
$$\sA_{E,\pi/\sM}\subset\sA_{E,\pi}, \quad
T_{\pi/\sM}\subset T_{X/\sM}\subset T_{\pi}$$ 
the pullback of $\pi^{-1}T_{S/\sM}\subset \pi^{-1}T_S$. Let
$$ ^{tr}\sA^{\bullet}_{E/\sM}:=(\,^{tr}\sA^{-1}_E\to\sA_{E,\pi/\sM}),\quad
\sA^{\bullet}_{E,\pi/\sM}:=(\sA_{E/S}\to\sA_{E,\pi/\sM}).$$
\end{notation}

\begin{prop}\label{prop4.3}
The exact sequences
\begin{gather}
\notag 0\to \omega_{X/S}[1]\to
\,^{tr}\sA^{\bullet}_E \to \sA^{\bullet}_{E, \pi}\to 0
\\ \notag
0\to \omega_{X/S}[1]\to
\,^{tr}\sA_{E/M}^{\bullet} \to \sA_{E, \pi/M}^{\bullet}\to 0
\end{gather}
have $0$-th direct images (via $\pi$) isomorphic
to
\begin{gather}
\notag 0\to \sO_S\to \sD^{\le 1}_S(\lambda_E)\to T_S\to 0
\\ \notag
0\to \sO_S\to \sD^{\le 1}_{S/M}(\lambda_E)\to T_{S/M}\to 0.
\end{gather}
\end{prop}

Furthermore, we need to review the
description of $^{tr}\sA^{\bullet}_E$ in terms 
of local coordinates, \cite{BS}, p. 660.    
Let $t$ be a local coordinate (along the fiber), 
and a trivilization $\sO^n_X\cong E$; $s$ 
a local coordinate on $S$.     
Note $t$ naturally induces local coordinates
$(t_1, t_2)$ around the diagonal of $X\times_S X$.   
One has isomorphisms
\begin{gather} \label{eq5f1}   
\sO_X\oplus {\rm Mat}_n(\sO_X)\oplus \sO_X\cong 
\,^{tr}\sA^{-1}_E;  
\\  
(\chi, B, \nu)\to \bigg[{\chi(t_1)\over (t_2-t_1)^2}+
{B(t_1)\over (t_2-t_1)}+ \nu(t_1)\bigg] dt_2;   
\notag 
\\
T_{\pi}\oplus {\rm Mat}_n(\sO_X)\cong\,  
^{tr}\sA^{0}_E=\sA_{E, \pi};  
(\tau, A)\to \tau(t,s)\partial_t+\mu(s)\partial_s+A. 
\notag    
\end{gather}
For different choices of coordinates and trivializations, 
there are formulas for transition functions,
namely the gauge change and coordinate change formulas, 
where $g\in {\rm GL}_r(\sO_X)$ and $y=y(t)$,
\begin{gather}\label{eq5f5}   
(\tau, A)
\buildrel{g}\over \longrightarrow 
(\tau, -\tau(g)g^{-1}+gAg^{-1}); 
\\
(\chi, B, \nu)
\buildrel{g}\over \longrightarrow 
\notag
\\
\bigg(\chi, -\chi g'g^{-1}+gBg^{-1}, 
{\rm Tr}\big(-{1\over 2}\chi g''g^{-1}+\chi(g'g^{-1})^2-Bg^{-1}g'\big)
+\nu\bigg); 
\notag
\\   
(\chi, B, \nu) \buildrel{y(t)}\over \longrightarrow 
\notag 
\\
\bigg(\chi {y'}^{-1}, B, r\chi \bigg({1\over 6}{y'''\over {y'}^{2}}
-{1\over 4}{{y''}^2\over {y'}^3}\bigg)+{1\over 2}{y''\over y'}{\rm Tr}B
+\nu y'\bigg).   
\notag
\end{gather}  

The main result of this section is
Theorem \ref{prop5.4} which
enables us to take care of Theorem \ref{prop0.2}.
Follow the notation of Section 2, let 
$$p: X=\sC\times_{\sM}S\to \sC$$
be the projection.  

\begin{defn}\label{propf6.1}
$\sE nd(E)^{-1}:={\rm res}^{-1}(\sE nd(E))
\subset\, ^{tr}\sA_E^{-1}$, and
$$\sG_E:=
{\rm res}^{-1}(\E)
\subset \sE nd(E)^{-1}$$
where $\sE nd(E)= \E \oplus \sO_X$
with its trace free part $\E:=\sE nd^0(E)$ and the trivial bundle $\sO_X$.
There are exact sequences
\begin{gather}
\notag
0\to\omega_{X/S}\to\sG_E\buildrel{\rm res}\over
\longrightarrow\E\to 0;
\\ \notag
0\to\omega_{X/S}\to\sE nd(E)^{-1}\buildrel{\rm res}\over \longrightarrow
\sE nd(E)\to 0.
\end{gather}
\end{defn}

Consider the natural morphism
$$S^2(\sG_E)\otimes_{\sO_X} T_{X/S}\xrightarrow{q}
S^2(\E)\otimes_{\sO_X}T_{X/S}\to 0$$
induced by $\sG_E \xrightarrow{\rm res}\E$.
Denote the kernel of $q$ by $\sK$.
There is a canonical isomorphism
$\iota:\sG_E\cong \sK$, that is 
$\iota(s)=\Sym^2(dt\otimes s)\otimes\partial_t$ locally.

\begin{defn}\label{defn6.2} 
$q^{-1}({\rm id}\otimes T_{X/S}):=S(\sG_E)$
where $\id\in S^2(\E)$ is the 
identity element. It follows that we have the exact sequence
\begin{gather}\label{eq6.1}0\to\sG_E\buildrel\iota\over
\longrightarrow S(\sG_E) \to T_{X/S}\to 0.
\end{gather}
Locally, for choosen coordinate and trivilization (cf. (\ref{eq5f1})),
any local section $s\in S(\sG_E)$ is of the form
$$s=\left(\aligned&\chi\sum_a(0,J_a, 0)\otimes (0,J_a,0)+
\sum_a\mu_a(0,J_a,0)\otimes(0,0,1)\\&+\sum_a\nu_a(0,0,1)\otimes(0,J_a,0)
+w(0,0,1)\otimes(0,0,1)\endaligned\right)\otimes\partial_t$$
where
$J_a$ is a (local) basis of $\E$ (which we assume 
to be orthonormal under $\rm Trace(\,\cdot\,)$)
such that $\sum_aJ_a\otimes J_a=\id$.
\end{defn}

We will need the Kodaira-Spencer maps
\begin{gather}\label{eq6.4.a}
{\rm KS}_{S}:
T_S\to R^1\pi_*\sA^0_{E/S}, \quad \sA^0_{E/S}=\sA_{E/S}/\sO_X;
\\ \notag
{\rm KS}_{\sM}: f^*T_{\sM}\to R^1\pi_*T_{X/S}
\end{gather}
They fit into the following commutative diagram
$$\CD
0@>>>R^1\pi_*(\E)@>>>R^1\pi_*\sA^0_{E/S}@>>>R^1\pi_*T_{X/S}@>>>0\\
@.@A{\rm KS}AA @A{\rm KS_S}AA @A{\rm KS_{\sM}}AA @.\\
0@>>>T_{S/\sM}@>>>T_S@>>>f^*T_{\sM}@>>>0
\endCD$$

\begin{remark}\label{prop6.4a1}
One way to see ${\rm KS}_{S}$ of \eqref{eq6.4.a} is
via the natural map (with cohomology) $\sA_{E,\pi}^{\bullet}\to \sA^0_{E/S}[1]$
from \eqref{eq4.8}; similarly ${\rm KS}_{\sM}$ is via
$(T_{\sC/\sM}\to T_{\pi})\to
T_{\sC/\sM}[1]$ (with cohomology) combined with its pullback via $f:S\to \sM$.
The diagram \eqref{eq6.4.a} via
natural maps $T_S\to f^*T_{\sM}$
and $R^1\pi_*\sA^0_{E/S}\to R^1\pi_*T_{X/S}$,
for Kodaira-Spencer maps commute.
\end{remark}

\begin{thm}\label{prop5.4}
i) If the Kodaira-Spencer map $KS:T_{S/\sM}\to R^1\pi_*(\E)$ is an
isomorphism, then there exists a canonical isomorphism
\begin{gather}
\notag \phi:R^1\pi_*(\sG_E)
\cong R^0\pi_*(\sE nd(\E)^{-1}\to \sA_{\E,\pi/\sM})
 \cong R^0\pi_*\,^{tr}\sA^{\bullet}_{\E/\sM}.
\end{gather}
ii) If ${\rm KS}_{S}$ is an isomorphism and shrink $\sM$ enough, 
then the above $\phi$ extends to an isomorphism
$$\phi=\phi_{-2r}:R^1\pi_*S(\sG_E) \cong 
R^0\pi_*(\,^{tr}\sA^{-1}_{\E}\to \sA_{\E,\pi}).$$
(The R.H.S. of i), ii) are canonically identified with
$\sD_{S/\sM}^{\le 1}(\lambda_{\E})$, $\sD_S^{\le 1}(\lambda_{\E})$ 
respectively,
cf. Proposition \ref{prop4.3}.)
\end{thm}

\begin{remark}\label{prop5.5}
$R^0\pi_*(\sE nd(\E)^{-1}\to \sA_{\E,\pi/\sM})
\cong R^0\pi_*(\sE nd(\E)^{-1}\to \sA_{\E,\pi})$
holds under $\sA_{\E,\pi/\sM} \to \sA_{\E,\pi}$, cf.
the proof of i) of Proposition \ref{prop6.4}.
\end{remark}

\begin{cor}\label{prop6.15}
Assumptions being as in ii) of Theorem \ref{prop5.4},
suppose $\lambda_{\E}\cong\Theta^{-\lambda}$ ($\lambda=2r$).
For $k\in \Z$ one has an isomorphism denoted by
$$\phi_{k}:
R^1\pi_*S(\sG_E)
\cong {\sD}^{\le 1}_S(\Theta^k), $$
extending Theorem \ref{prop5.4}.
(If $k=-\lambda$, write $\phi$ for $\phi_{-\lambda}$.)
\end{cor}

\begin{prop}\label{prop5.6}
i) The morphism $\E\to \sE nd^0(\E)$
induced by the adjoint representation extends naturally
to 
a canonical morphism
${\rm ad}:\sA_E \to \sA_{\E}$ (preserving algebra structures).
ii) The morphism ${\rm ad}$ has a natural lifting
$\ov{\rm ad}:\sG_E\to \sG_{\E}$, which induces $(2r)\cdot{\rm id}$
on $\omega_{X/S}$. 
\end{prop}
\begin{proof}  i) For any $D\in \sA_E$, $L\in\E$,
${\rm ad}(D)(L):=D\circ L-L\circ D$ is a section of $\E$ (note that
${\rm Tr}(D\circ L-L\circ D)=\e_1(D){\rm Tr}(L)$). Thus
${\rm ad}(D)$ defines a map $\E\to\E$, which is a differential
operator since ${\rm ad}(D)(\lambda\cdot L)=\lambda\cdot {\rm ad}(D)(L)+
\e_1(D)(\lambda)\cdot L$. 

ii)
This can be proved
via the local formulas given in (\ref{eq5f1})
Namely, a local element of $\sG_E$ is expressed
as $(0,B,\nu)$ (with $B\in {\rm Mat}_r({\sO_X})$
and $\nu\in \sO_X$).  Define a
lifting by sending
$(0,B, \nu)$ to $(0,{\rm ad}B, 2r\nu)$.
We will show by using formulas in (\ref{eq5f5}) that
the above lifting is in fact globally defined. 
Using $\chi=0$, ${\rm Tr}B=0$
and ${\rm Tr}({\rm ad}B)=0$ in (\ref{eq5f5}), we have
\begin{gather}\label{eq5f6}   
(0, B, \nu)
\buildrel{g}\over \longrightarrow 
\big(0, gBg^{-1}, {\rm Tr}(-Bg^{-1}g')+\nu\big); 
\notag
\\   
(0, B, \nu) \buildrel{y(t)}\over \longrightarrow 
\big(0, B, \nu y'\big);    
\notag
\\
(0, {\rm ad}B, \nu)
\buildrel{g}\over \longrightarrow 
\Big(0, {\rm ad}(gBg^{-1}), {\rm Tr}\big(-{\rm ad}B\,{\rm ad}(g^{-1}g')\big)
+\nu\Big); 
\notag
\\ 
(0, {\rm ad}B, \nu)
\buildrel{y(t)}\over \longrightarrow 
(0, {\rm ad}B, \nu y'). 
\notag
\end{gather}
If change the trivilization of $E$ by $g$, the induced trivilization of
$\sE nd(E)$ will be changed by $e_g:\M_r(\sO_X)\to\M_r(\sO_X)$, where
$e_g(B)=gBg^{-1}$. It is easy to check that 
$e_g^{-1}e_g'={\rm ad}(g^{-1}g')$, thus we obtain    
 the term ${\rm Tr}\big(-{\rm ad}B\,{\rm ad}(g^{-1}g')\big)$
in the 3rd row above.  
One knows that 
\begin{gather}\label{eq5f18}
{\rm Tr}({\rm ad}M\,{\rm ad}N)= 
2r{\rm Tr}(MN)
\end{gather}
for traceless matrices $M$, $N$ of rank $r$.      
Let $(g^{-1}g')_0$ be the traceless compoment of $g^{-1}g'$.   
Note $\Tr B=0$, ${\rm ad}(g^{-1}g')={\rm ad}((g^{-1}g')_0)$,   
and ${\rm Tr}(Bg^{-1}g')={\rm Tr}(B(g^{-1}g')_0)$ etc.  
It follows that the morphism $(0,B, \nu)\to (0, {\rm ad}B, 2r\nu)$ 
as given is well-defined (globally).  
\end{proof}

\begin{lem}\label{prop5.7}
$R^0\pi_*\sA_{E,\pi}\cong \sO_S$ and $R^0\pi_*\sA^0_{E,\pi}=0$
provided ${\rm KS}_{S}$ being injective, where 
$\sA^0_{E,\pi}:=\sA_{E,\pi}/\sO_X$.
\end{lem}
\begin{proof}
It suffices to prove,
by using an exact sequence similar to
\eqref{eq4.2} (with $T_{\bullet}$
replaced by $\sA_{\bullet}$),
that { i)} $\pi_*\pi^{-1}T_{S}\to R^1\pi_*\sA_{E/S}$
is injective and { ii)} $\pi_*\sA_{E/S}\cong \sO_S$.
But the map in { i)} composed with
$R^1\pi_*\sA_{E/S}\to R^1\pi_*\sA^0_{E/S}$
is nothing but ${\rm KS}_{S}$,
hence { i)}.   { ii)} is
from $\pi_*T_{X/S}=0$
(genus $\ge 2$) and $\pi_*\sE nd(E)=\sO_S$
($E$ is fiberwise stable).
\end{proof}

We are ready to give a proof of i) of Theorem \ref{prop5.4}.
\begin{proof}
Firstly, we remark that morphisms in 
Proposition \ref{prop5.6} make the following
diagram of complexes being commutative 
$$\CD
{(\,^{tr}\sA^{-1}_{\E}\buildrel{\rm res}\over\longrightarrow\sA_{\E,\pi/\sM})}
@>>>{(\sA_{\E/S}\to\sA_{\E,\pi/\sM})}\\
@AAA  @AAA \\
{(\sG_E\buildrel{\rm res}\over\longrightarrow\sA_{E,\pi/\sM})}@>>>
{(\sA^0_{E/S}\to\sA^0_{E,\pi/\sM})}
\endCD$$
Secondly, we observe that the commutative diagram
$$\CD
0@>>>\E@>>>\sA^0_{E,\pi/\sM}@>>>T_{\pi/\sM}@>>>0\\
@.@| @AAA @AAA @.\\
0@>>>\E@>>>\sA^0_{E/S}@>>> T_{X/S}@>>>0
\endCD$$
induces a commutative diagram
$$\CD
T_{S/\sM}@=\pi_*T_{\pi/\sM}@>{\rm KS}>>R^1\pi_*(\E)@>>>R^1\pi_*\sA^0_{E,\pi/\sM}\\
@.@AAA @| @AAA \\
@.0@>>>R^1\pi_*(\E)@>>>R^1\pi_*\sA^0_{E/S}
\endCD$$
Thus $R^1\pi_*(\E)$ vanishes in $R^1\pi_*\sA^0_{E,\pi/\sM}$ since KS is an
isomorphism.

We construct $\phi$ for any affine open set $U^i\subset S$. 
Let $\{U_i,\,\dot X_i\}$ be an affine covering of $\pi^{-1}(U^i)$, let 
$\dot U_i=U_i\cap\dot X_i$. 
For any $\check {\rm C}$ech cocycle
$r_{\dot U_i}\in \sG_E(\dot U_i)$
in $C^1(\sG_E)$, the class $[{\rm res}(r_{\dot U_i})]\in R^1\pi_*(\E)(U^i)$ 
vanishes in $R^1\pi_*\sA^0_{E,\pi/\sM}(U^i)$. Thus there exists
$\tau_{\dot X_i}\in \sA^0_{E,\pi/\sM}(\dot X_i)$,
$\tau_{U_i}\in  \sA^0_{E,\pi/\sM}(U_i)$ 
such that
$\tau_{\dot X_i}- \tau_{U_i}={\rm res}\,r_{\dot U_i}$ on $\dot U_i$. 
For given $r_{\dot U_i}$, the choice of $\tau_{\dot X_i}$ and $\tau_{U_i}$ 
is unique (by Lemma \ref{prop5.7}). Then
\begin{gather}\label{eq5.2}
\{{\rm ad}(\tau_{\dot X_i}),{\rm ad}(\tau_{U_i});\ov{\rm ad}(r_{\dot U_i})\}
\end{gather}
is a $\check {\rm C}$ech cocycle in 
$ C^0(\sG_{\E}\to \sA_{\E,\pi/\sM})$ (cf. \cite{BS}, p. 673).
It is easily checked that the assignment
\begin{gather}\label{eq5.4.3}\tilde\phi:r_{\dot U_i}\to
\{{\rm ad}(\tau_{\dot X_i}), {\rm ad}(\tau_{U_i}); \ov{\rm ad}(r_{\dot U_i})\}
\quad \in 
C^0(\,^{tr}\sA^{\bullet}_{\E/\sM})\end{gather}
preserves the respective coboundaries. Hence it descends to a map 
$$\phi: R^1\pi_*(\sG_E) \to R^0\pi_*\,^{tr}\sA^{\bullet}_{\E/\sM}.$$
By the same way, we construct 
$\vartheta:R^1\pi_*(\E)\to R^0\pi_*\sA^{\bullet}_{\E,\pi/\sM}$ such
that 
$$\CD
0@>>>R^1\pi_*\omega_{X/S}@>>>
R^1\pi_*(\sG_E)@>{\rm res}>>
R^1\pi_*(\E)@>>>0\\
@. @V(2r)\cdot{\rm id}VV     @V\phi VV @V\vartheta VV  @.  \\
0@>>>R^0\pi_*\omega_{X/S}[1]@>>>
R^0\pi_*\,^{tr}\sA^{\bullet}_{\E/\sM}
@>{\rm res}>>R^0\pi_*\sA^{\bullet}_{\E,\pi/\sM}@>>>0
\endCD$$
is commutative. The map $\vartheta$ is the composition of  
$${\rm KS}^{-1}: R^1\pi_*(\E)\to T_{S/\sM}=
R^0\pi_*(\sA^0_{E/S}\to\sA^0_{E,\pi/\sM})$$
and the map 
$R^0\pi_*(\sA^0_{E/S}\to\sA^0_{E,\pi/\sM})\to 
R^0\pi_*\sA^{\bullet}_{\E,\pi/\sM}$, which is induced by 
the qusi-isomorphism of complexes at the begining of our proof. Thus
$\vartheta$ is an isomorphism, then $\phi$ has to be an isomorphism.
\end{proof}

Both $R^1\pi_*(\sG_E)$ and
$R^0\pi_*(\sE nd(\E)^{-1}\to \sA_{\E,\pi/\sM})
=R^0\pi_*\,^{tr}\sA^{\bullet}_{\E/\sM}$ define
extension classes in $H^1(S,\Omega_{S/\sM})$.
One has
(by the preceding proof combined with
proof of Proposition \ref{prop5.6} for 
the constant $2r$) 

\begin{cor}\label{prop5.9}
The extension classes (${\rm e.c.} [\bullet]$ for short) 
satisfy $${\rm e.c.}[R^1\pi_*(\sG_E)]=
{1\over 2r} {\rm e.c.}[\sD^{\le 1}_{S/\sM}(\lambda_{\E})]$$
in $H^1(S, \Omega_{S/\sM})$, where  
$\sD^{\le 1}_{S/\sM}(\lambda_{\E})= 
R^0\pi_*(\sE nd(\E)^{-1}\to \sA_{\E,\pi/\sM})$.  
\end{cor}

For ii) of the theorem, our proof 
will need the following result.

\begin{prop}\label{props6.12} There exisits 
 $\ov{\rm res}: S(\sG_E)\to \sA^0_{E/S}$ such that
$$\CD
0@>>>\omega_{X/S}@>>> S(\sG_E)@>\ov{\rm res}>>\sA^0_{E/S}@>>>0\\
@. @VVV @| @VVV @.\\
0@>>>\sG_E@>\iota>> S(\sG_E)@>>>T_{X/S}@>>>0
\endCD$$
is commutative.  
If ${\rm KS}_S$ in \eqref{eq6.4.a} is an isomorphism,
then 
\begin{gather} \notag 0\to \sO_S\to
R^1\pi_* S(\sG_E) \to T_S\to 0.
\end{gather}
\end{prop}
\begin{proof} For any local section $s\in S(\sG_E)$ in 
Definition \ref{defn6.2}, we define
$$\ov{\rm res}(s)=(\chi\partial_t, \frac{1}{2}\sum_a(\mu_a+\nu_a)J_a),$$
which is independent of the choice of $\{J_a\}$, thus well-defined locally.
To show it well-defined globally, we need to check the invarience under
gauge and coordinate changes. The invarience under local coordinate changes
is straightforward. Under the gauge change $g\in {\rm GL}_r(\sO_X)$,
the section $s$ becomes into $s_g\otimes\partial_t$, where $s_g$ is
$$\aligned &\chi\sum_a(0,gJ_ag^{-1},0)^{\otimes 2}+
\sum_a(\mu_a-\chi{\rm Tr}(J_ag^{-1}g'))(0,gJ_ag^{-1},0)\otimes (0,0,1)\\
&+\sum_a(\nu_a-\chi{\rm Tr}(J_ag^{-1}g'))(0,0,1)\otimes(0,gJ_ag^{-1},0)\\
&+\left(\chi\sum_a{\rm Tr}(J_ag^{-1}g')^2-\sum_a(\mu_a+\nu_a)
{\rm Tr}(J_ag^{-1}g')+w\right) (0,0,1)^{\otimes 2}\endaligned.$$
Then $\ov{\rm res}(s_g\otimes\partial_t)=
(\chi\partial_t,\frac{1}{2}\sum_a
(\mu_a+\nu_a-2\chi{\rm Tr}(J_ag^{-1}g'))gJ_ag^{-1})$ coincides with
$(\chi\partial_t,-\chi\partial_t(g)g^{-1}+\frac{1}{2}
\sum_a(\mu_a+\nu_a)gJ_ag^{-1})$ in $\sA^0_{E/S}=\sA_{E/S}/\sO_X$ since
$\chi\partial_t(g)g^{-1}=\chi g'g^{-1}=
\sum_a\chi{\rm Tr}(J_ag^{-1}g'))gJ_ag^{-1}$ modulo $\sO_X$.
Thus $\ov{\rm res}$ is gauge invarient and defined globally.
Then the rest of this proposition is obvious.
\end{proof}

\begin{remark}\label{propsym}
i) $\Sym^2(a\otimes b)={1\over 2} (a\otimes b+b\otimes a)$
for $a\,, b\in \sG_E$.
ii) Using Proposition \ref{props6.12}  
and assuming $R^1f_*\sO_S=0$ 
one has a quick interpretation of 
ii) of Theorem \ref{prop5.4}.      
Note both $R^1\pi_* S(\sG_E)$
and $R^0\pi_*\,^{tr}\sA^{\bullet}_{\E}\cong \sD^{\le 1}_S(\lambda_{\E})$
contain a subsheaf $R^1\pi_*(\sG_E)\cong
\sD^{\le 1}_{S/\sM}(\lambda_{\E})$.
Given two extensions $\sF$, $\sF'$
of $T_S$ by $\sO_S$ suppose their subsheafs
with symbolic part in $T_{S/\sM}$ are isomorphic,
then $\sF\cong \sF'$ (non-canonically) provided
$R^1f_*\sO_S=0$ (since $H^1(S, \Omega_S)\to H^1(S, \Omega_{S/\sM})$
is injective if $R^1f_*\sO_S=0$ (and $\sM$ is affine)).
The ii) of Theorem \ref{prop5.4}
just proves an isomorphism of this kind without
reference to $R^1f_*\sO_S=0$.
\end{remark}

In what follows, for simplicity, we will cover $\sC$ by two affine open sets
$V$ and $\dot\sC$ (shrink $\sM$ if necessary). Then we choose
and fix the local coordinates (along the fibre) on $V$ and $\dot\sC$.
Let $U=p^{-1}(V)$, $\dot X=p^{-1}(\dot\sC)$ and $\{U^i\}_i$ be an affine
covering of $S$. Let $U_i:=U\cap\pi^{-1}(U^i)$ and 
$\dot X_i=\dot X\cap\pi^{-1}(U^i).$ Then it is important that on each $U_i$
(resp. $\dot X_i$) we use the local coordinate pulling back from $V$ 
(resp. $\dot\sC$). Thus any construction on $U$ (resp. $\dot X$)
using the local description (\ref{eq5f1}) only depends on the trivialization
of $E$ over $U_i$ (resp. $\dot X_i$). We start with the
construction of 
$\gamma_{\E,U}:S(\sG_E)|_U \to \,^{tr}\sA^{-1}_{\E}|_U$ (resp.
$\gamma_{\E,\dot X}:S(\sG_E)|_{\dot X} \to \,^{tr}\sA^{-1}_{\E}|_{\dot X}$)
such that the following diagram $(*)$ is commutative over $U$ 
(resp. over $\dot X$) 
$$\CD
0@>>> \omega_{X/S}@>>>\,^{tr}\sA^{-1}_{\E}@>{\rm res}>>\sA_{\E/S}@>>>0\\
@.  @A r\cdot{\rm id}AA @A\gamma_{\E,U} AA @A{\rm ad} AA @.\\
0@>>>\omega_{X/S}@>>> S(\sG_E)@>\ov{\rm res}>>\sA^0_{E/S}@>>>0\\
@. @VVV @| @VVV @.\\
0@>>>\sG_E@>\iota>> S(\sG_E)@>>>T_{X/S}@>>>0
\endCD$$
where ${\rm ad}:\sA^0_{E/S}\to\sA_{\E/S}$ is induced by the morphism in 
Proposition \ref{prop5.6} i) that maps $\sO_X$ to zero. Note that,
except $\gamma_{\E,U}$ (resp. $\gamma_{\E,\dot X}$), other morphisms
in the diagram are well defined over the global $\sM$ (i.e. need not to
shrink $\sM$).
 
On each $U_i$, fix a trivilization
of $E$ on $U_i$ and use the pullback coordinate of $V$, 
we define the morphism by using the local description \eqref{eq5f1}. 
For any local section 
$$\alpha=\left(\aligned&\chi\sum_a(0,J_a, 0)\otimes (0,J_a,0)+
\sum_a\mu_a(0,J_a,0)\otimes(0,0,1)\\&+\sum_a\nu_a(0,0,1)\otimes(0,J_a,0)
+w(0,0,1)\otimes(0,0,1)\endaligned\right)\otimes\partial_t$$
one defines that (with $r$ the rank of $E$)
\begin{gather}
\label{eq5.4.10}
\gamma_{\E,U}(\alpha)=(\chi, {1\over 2} \sum_a(\mu_a+\nu_a){\rm ad}J_a,
rw)\in\, ^{tr}\sA_{\E}^{-1}(U_i).
\end{gather}

\begin{lem}\label{prop5.20}
The assignment $\alpha\to \gamma_{\E,U}(\alpha)$ 
constructed above is gauge-invariant 
(it is not, however, independent of the choice of coordinate on $V$).     
Equivalently, under another choice of 
trivialization of $E$ (on $U_i$), such that $J_a\to gJ_ag^{-1}$ 
and $\alpha\to \alpha_g$, 
the assignment remains unchanged, i.e. $\gamma_{\E,U}(\alpha_g)$ 
is obtained as the $e_g$-transformation of $\gamma_{\E,U}(\alpha)$,
where $e_g$ is the gauge of $\E$ induced by $g$ 
(cf. Proposition \ref{prop5.6}).   
\end{lem}
\begin{proof} As in the proof of Proposition \ref{props6.12},
$\alpha_g=s_g\otimes\partial_t$. Then
$$\aligned\gamma_{\E,U}(\alpha_g)=(&\chi,\frac{1}{2}
\sum_a(\mu_a+\nu_a-2\chi{\rm Tr}(J_ag^{-1}g')){\rm ad}(gJ_ag^{-1}),\\
&r\chi\sum_a{\rm Tr}(J_ag^{-1}g')^2-r\sum_a(\mu_a+\nu_a)
{\rm Tr}(J_ag^{-1}g')+rw).\endaligned$$
The $e_g$-transformation $\gamma_{\E,U}(\alpha)^g$ of 
$\gamma_{\E,U}(\alpha)$ is
$$\left(\aligned&\chi,-\chi e_g'e_g^{-1}+\frac{1}{2}
\sum_a(\mu_a+\nu_a)e_g{\rm ad}(J_a)e_g^{-1},\quad rw+\\&
{\rm Tr}(-\frac{1}{2}\chi e_g''e_g^{-1}+\chi (e_g'e_g^{-1})^2-
\frac{1}{2}\sum_a(\mu_a+\nu_a){\rm ad}(J_a)e_g^{-1}e_g')\endaligned\right)$$
Recall that $e_g:\M_r(\sO_X)\to\M_r(\sO_X)$ means $e_g(B)=gBg^{-1}$,
we have $e_g{\rm ad}(J_a)e_g^{-1}={\rm ad}(gJ_ag^{-1})$. Thus the second
components of $\gamma_{\E,U}(\alpha_g)$ and $\gamma_{\E,U}(\alpha)^g$ 
will coincide if 
$e_g'e_g^{-1}=\sum_a{\rm Tr}(J_ag^{-1}g'){\rm ad}(gJ_ag^{-1}),$
which is true since $e_g'e_g^{-1}={\rm ad}(g'g^{-1})$,  
$e_g^{-1}e_g'={\rm ad}(g^{-1}g')$. To finish the proof, we will show that
their third components coincide. Since
$$\frac{1}{2}\sum_a(\mu_a+\nu_a){\rm Tr}({\rm ad}(J_a)e_g^{-1}e_g')
=r\sum_a(\mu_a+\nu_a){\rm Tr}(J_ag^{-1}g'),$$
it will be done if one can show the following identity
\begin{gather}\label{eq5f17}
r\sum_a{\rm Tr}(J_ag^{-1}g')^2=
{\rm Tr}((e_g'e_g^{-1})^2-\frac{1}{2}e_g''e_g^{-1}).
\end{gather}
Write $e_g''e_g^{-1}=(e_g'e_g^{-1})'-e_g'(e_g^{-1})'
=(e_g'e_g^{-1})'+ e_g'e_g^{-1}e_g'e_g^{-1}$, then
$${\rm Tr}(e_g''e_g^{-1})={\rm Tr}((e_g'e_g^{-1})^2)=
{\rm Tr}({\rm ad}(g'g^{-1}){\rm ad}(g'g^{-1}))$$
(using $(e_ge_g^{-1})'=0$ and 
${\rm Tr}((e_g'e_g^{-1})')={\rm Tr}(e_g'e_g^{-1})'=0$ here). Let
$(g'g^{-1})_0$ be the traceless part
of $g'g^{-1}$. Then R.H.S of (\ref{eq5f17}) equals to 
$$\frac{1}{2}{\rm Tr}({\rm ad}(g'g^{-1}){\rm ad}(g'g^{-1}))=
r{\rm Tr}((g'g^{-1})_0(g'g^{-1})_0).$$
By the choice of $\{J_a\}_a$, we have 
$(g^{-1}g')_0=\sum_a{\Tr}(J_ag^{-1}g')J_a$. Then 
L.H.S of (\ref{eq5f17}) equals to
$r{\rm Tr}((g^{-1}g')_0g^{-1}g')=r{\rm Tr}((g^{-1}g')_0(g^{-1}g')_0).$
Thus (\ref{eq5f17}) is true since $g(g^{-1}g')_0g^{-1}=(g'g^{-1})_0$.
We are done.
\end{proof}

We have constructed $\gamma_{\E,U}$ (the 
construction of $\gamma_{\E,\dot X}$ is similar) such that 
the above diagram $(*)$
is commutative over $U$ (resp. over $\dot X$). It is also easy to see
that $\gamma_{\E,U}$, $\gamma_{\E,\dot X}$ induce (through $\iota$) 
$\frac{1}{2}\ov{\rm ad}$ on
$\sG_E$. However, 
$\gamma_{\dot U}:=\gamma_{\E,U}-\gamma_{\E,\dot X}$ may not vanish on 
$S(\sG_E)|_{\dot U}$ (but vanishes on $\sG_E|_{\dot U}$),
which defines a morphism
$\gamma_{\dot U}:S(\sG_E)|_{\dot U}\to \omega_{X/S}|_{\dot U}$
that induces $\ov\gamma_{\dot U}:T_{X/S}|_{\dot U}\to\omega_{X/S}|_{\dot U}$,
i.e., a section 
$$\ov \gamma_{\dot U}\in \sH om(T_{X/S},\omega_{X/S})(\dot U)=
\omega_{X/S}^2(p^{-1}(\dot V))$$
where $\dot V:=V\cap\dot\sC$. $[\ov\gamma_{\dot U}]$ defines a 
class of $H^1(\sC,p_*\omega^2_{X/S})=H^1(\sC,\omega^2_{\sC/\sM})$, which
vanishes since we assume $\sM$ being affine. Thus there exist
$$\ov\psi_U\in \sH om(T_{X/S},\omega_{X/S})(U),\quad  
\ov\psi_{\dot X}\in \sH om(T_{X/S},\omega_{X/S})(\dot X)$$
such that $\ov \gamma_{\dot U}=\ov\psi_U-\ov\psi_{\dot X}$. Let $\psi_U$,
$\psi_{\dot X}$ denote the induced morphisms
$$\psi_U: S(\sG_E)|_U\to S(\sG_E)/\sG_E|_U\cong T_{X/S}|_U
\buildrel{\ov\psi_U}\over\longrightarrow \omega_{X/S}|_U$$ 
$$\psi_{\dot X}: S(\sG_E)|_{\dot X}\to S(\sG_E)/\sG_E|_{\dot X}
\cong T_{X/S}|_{\dot X}
\buildrel{\ov\psi_{\dot X}}\over\longrightarrow \omega_{X/S}|_{\dot X}.$$
Then it is easy to see that 
$\gamma_{\dot U}=\gamma_{\E,U}-\gamma_{\E,\dot X}=\psi_U-\psi_{\dot X}$. 
Let $$\beta_U:=\gamma_{\E,U}-\psi_U,\quad \beta_{\dot X}:=
\gamma_{\E,\dot X}-\psi_{\dot X}.$$
Thus, by shrinking $\sM$, we have proved the following

\begin{prop}\label{prop5.14}
The $\beta_U$ and $\beta_{\dot X}$ define a morphism
$$\beta:S(\sG_E)\to \,^{tr}\sA^{-1}_{\E},$$
which induces (through $\iota$) $\frac{1}{2}\ov{\rm ad}$ on $\sG_E$, 
such that the following diagram 
is commutative
$$\CD
0@>>> \omega_{X/S}@>>>\,^{tr}\sA^{-1}_{\E}@>{\rm res}>>\sA_{\E/S}@>>>0\\
@.  @A r\cdot{\rm id}AA @A\beta AA @A{\rm ad} AA @.\\
0@>>>\omega_{X/S}@>>> S(\sG_E)@>\ov{\rm res}>>\sA^0_{E/S}@>>>0\\
@. @VVV @| @VVV @.\\
0@>>>\sG_E@>\iota>> S(\sG_E)@>>>T_{X/S}@>>>0
\endCD$$
\end{prop}

We shall now prove ii) of Theorem \ref{prop5.4}.

\begin{proof}
The proof is similar to that of i) of the theorem. By 
\begin{gather}\label{eq5.4.6} 
0\to \sA^0_{E/S}\to \sA^0_{E,\pi}\to \pi^{-1}T_S\to 0  
\end{gather} 
and by identifying its connecting map $T_S\to R^1\pi_*\sA^0_{E/S}$
with the Kodaira-Spencer map $KS_S$ (which will be treated
more generally for $\log D$ in Proposition \ref{prop6.4}),
we see that $R^1\pi_*\sA^0_{E/S}$ vanishes in $R^1\pi_*\sA^0_{E,\pi}$
if $KS_S$ is an isomorphism.
Thus, for any $\alpha\in S(\sG_E)(\dot U_i)$, there exist
$\tau_{\dot X_i}\in \sA^0_{E,\pi}(\dot X_i)$ and 
$\tau_{U_i}\in \sA^0_{E,\pi}(U_i)$ such that
$\tau_{\dot X_i}-\tau_{U_i}=\ov{\rm res}(\alpha).$
Then, by Proposition \ref{prop5.14}, we see that
${\rm res}(\beta(\alpha))={\rm ad}(\ov{\rm res}(\alpha)).$
Thus 
$$\ov\phi(\alpha):=
\{{\rm ad}(\tau_{\dot X_i}),{\rm ad}(\tau_{U_i}); \beta(\alpha)\}$$
is a cocycle in 
$C^0(\,^{tr}\sA^{-1}_{\E}\to \sA^0_{\E,\pi})$. It is clear that $\ov\phi$ 
induces a morphism
$$\phi: R^1\pi_*S(\sG_E) \to 
R^0\pi_*(\,^{tr}\sA^{-1}_{\E}\to \sA^0_{\E,\pi})=
R^0\pi_*\,^{tr}\sA^{\bullet}_{\E}.$$
Similarly, we can construct 
$$\vartheta: R^1\pi_*\sA^0_{E/S}\to R^0\pi_*(\sA^0_{E/S}\to\sA^0_{E,\pi})
\to R^0\pi_*\,^{tr}\sA^{\bullet}_{\E},$$
which is an isomorphism such that the following diagram  
$$\CD
0@>>>R^1\pi_*\omega_{X/S}@>>>R^1\pi_*S(\sG_E)@>\ov{\rm res}>>
R^1\pi_*\sA^0_{E/S}@>>>0\\
@. @V r\cdot{\rm id}VV     @V\phi VV @V\vartheta VV  @.  \\
0@>>>R^0\pi_*\omega_{X/S}[1]@>>>
R^0\pi_*\,^{tr}\sA^{\bullet}_{\E}
@>{\rm res}>>R^0\pi_*\sA^{\bullet}_{\E,\pi}@>>>0
\endCD$$
is commutaive. Thus $\phi$ must be an isomorphism.
\end{proof}

\begin{remark}\label{rem3.18} 
The $\phi$ in Theorem \ref{prop5.4} i) is defined globally over
$\sM$, but the one in Theorem \ref{prop5.4} ii) is defined only over
an open set of $\sM$. More precisely, there is an affine covering
$\{\sM_i\}_{i\in I}$ of $\sM$ such that on each 
$\sC\times_{\sM}f^{-1}(\sM_i)$ we can choose a 
$\beta_i:S(\sG_E)\to \,^{tr}\sA^{-1}_{\E}$
as the $\beta$ in Proposition \ref{prop5.14}. Then, by using $\beta_i$ and
Theorem \ref{prop5.4} ii), we get the isomorphism 
$\phi_i: R^1\pi_*S(\sG_E)\to
R^0\pi_*\,^{tr}\sA^{\bullet}_{\E}=
\sD_S^{\le 1}(\lambda_{\E})$
on $f^{-1}(\sM_i)$.
For another choice $\{\beta_i'\}_{i\in I}$, the map 
$\beta_i-\beta_i': S(\sG_E)\to\omega_{X/S}$ induces
$$\phi_i-\phi_i':R^1\pi_*S(\sG_E) \to R^0\pi_*\omega_{X/S}[1]=\sO_S$$
on $f^{-1}(\sM_i)$, which vanishes on $R^1\pi_*(\sG_E)$, thus 
$(\phi_i-\phi_i')\in \Omega^1_{\sM}(\sM_i)$.
Similarly, on $f^{-1}(\sM_i\cap\sM_j)$, $\phi_{ij}:=\phi_i-\phi_j$ induces
$\ov\phi_{ij}\in \Omega^1_{\sM}(\sM_i\cap\sM_j).$ Thus $\{\ov\phi_{ij}\}$
defines a class in ${\rm H}^1(\sM,\Omega^1_{\sM})$. 
\end{remark}

To conclude this section, we describe the connecting maps
$\tilde\delta$ and prove Lemma \ref{prop0} (cf. Lemma \ref{prop5.13}).

\begin{lem}\label{prop5.13} The map $\tilde\delta$ induces the identity
map on $R^1\pi_*(\sG_E)$. More precisely, 
$\tilde\delta\circ\iota(\bullet)=\Sym^2((\bullet)\otimes 1).$
\end{lem}
\begin{proof} It is known (see Proposition 4.2 of \cite{R}) that the
connecting map 
$$\tilde\delta: R^1\pi_*(\sG_E\boxtimes\sG_E(\Delta)|_{\Delta})\to
R^2(\pi\times\pi)_*(\sG_E\boxtimes\sG_E)$$
is dual (under Serre duality) to the restriction map 
$$r: (\pi\times\pi)_*(\sG_E^*\boxtimes\sG_E^*\otimes
\omega_{X/S}\boxtimes\omega_{X/S})\to \pi_*(\sG_E^*\otimes\sG_E^*\otimes
\omega_{\Delta/S}^2).$$
Work locally on $S$, we can assume that $X$ is covered by affine open
sets $U$ and $\dot X$. Let $w\in \omega_{X/S}(\dot U)$ be a base of $\omega_{X/S}$
on $\dot U:=U\cap\dot X$ and $w^*\in T_{X/S}(\dot U)$ be its dual base.
Then, for any $\alpha\in \sG_E(\dot U)$, 
$$\iota([\alpha])=[\Sym^2(\alpha\otimes w)\otimes w^*]
=\frac{1}{2}[(\alpha\otimes w+w\otimes\alpha)\otimes w^*].$$
We use the following identification
$$(\pi\times\pi)_*(\sG_E^*\boxtimes\sG_E^*\otimes
\omega_{X/S}\boxtimes\omega_{X/S})\cong 
\pi_*(\sG_E^*\otimes\omega_{X/S})\otimes\pi_*(\sG_E^*\otimes\omega_{X/S}).$$ 
For any $\beta_i\in\pi_*(\sG_E^*\otimes\omega_{X/S})$ ($i=1,2$), let
$\beta_i|_{\dot U}=s_i\otimes w$, where $s_i\in\sG_E^*(\dot U)$ ($i=1,2$).
Then $r(\beta_1\otimes\beta_2)|_{\dot U}=s_1\otimes s_2\otimes w\otimes w$. 
Thus   
\begin{gather}\label{eq5.131}
<\tilde\delta(\iota([\alpha])),\beta_1\otimes\beta_2>
=\frac{1}{2}[(s_1(\alpha)s_2(w)+s_1(w)s_2(\alpha))\cdot w].
\end{gather}

By $R^1\pi_*\omega_{X/S}=\sO_S$, let 
$[f\cdot w]=1\in\sO_S$ ($f\in\sO_X(\dot U)$), we have
\begin{gather}\label{eq5.132} 
<\Sym^2([\alpha]\otimes [f\cdot w]),\beta_1\otimes\beta_2>
\\
\notag
=\frac{1}{2}([s_1(\alpha)\cdot w][fs_2(w)\cdot w]+
[fs_1(w)\cdot w][s_2(\alpha)\cdot w]).
\end{gather}
Note that $s_i(w)=\beta_i|_{\dot U}(w\otimes w^*)$ ($i=1,2$), we can see
that $s_i(w)=\beta_i({\rm id})|_{\dot U}$, where the (global) section
${\rm id}$ is the image of $1$ under $\sO_X \to \sG_E\otimes\omega_{X/S}^*$.
Thus $s_i(w)\in\sO_S$ since $\beta_i$ ($i=1,2$) are global sections of
$(\sG_E\otimes\omega^*_{X/S})^*$. Then 
$[s_1(\alpha)s_2(w)\cdot w]=[s_1(\alpha)\cdot w][fs_2(w)\cdot w]$ and
$[s_2(\alpha)s_1(w)\cdot w]=[s_2(\alpha)\cdot w][fs_1(w)\cdot w]$, which means that
$$<\tilde\delta(\iota([\alpha])),\bullet>=
<\Sym^2([\alpha]\otimes 1),\bullet>.$$
\end{proof}

\section{The Generalization to Singular Cases}

Let $\pi:\ov X\to \ov S$ be a proper morphism of relative
dimension $1$ between smooth varieties in
characteristic 0 such that each fiber has at most
ordinary double points as singularities.  Let
$\ov f:\ov S\to \ov\sM$ be a smooth morphism where $\ov\sM$ is a smooth
variety. Let $B=\ov\sM\setminus\sM$ and $W={\ov f}^{-1}(B)$ such that
$D=\pi^{-1}(W)$ consists precisely of singular fibres. 
As before let $\omega_{\ov X/\ov S}$ be the relative dualizing
sheaf (which is locally free as is well known).
Let $T_{\ov\sM}(B)\subset T_{\ov\sM}$ be the subalgebra of
vector fields that preserve $B$ (cf. \cite{BS}, Section 6).

\begin{notation}
In the notation of Section 3, define the following
\begin{gather}\notag
T_{\ov f}(\log D)=d({\ov f}\circ\pi)^{-1}
((\ov{f}\circ\pi)^{-1}T_{\ov\sM}(B))\subset T_{\ov X};\\
\notag T_{\pi}(\log D)=T_{\ov f}(\log D)\cap d\pi^{-1}
(\pi^{-1}T_{\ov S})\subset T_{\ov X};
\\ \notag
T_{\ov S}(\log W)=d\pi (T_{\pi}(\log D))\subset T_{\ov S};\\
\notag T_{\ov X/\ov S}(\log D)=T_{\ov X/\ov S}\cap T_{\pi}(\log D).
\end{gather}
\end{notation}
\begin{notation}
Let $E$ be a vector bundle on $\ov X$.
Define the following
\begin{gather}\notag
\sA_{E,\pi}(\log D) =\epsilon^{-1}T_{\pi}(\log D)\subset \sA_E;\\
\notag\sA_{E/\ov S}(\log D) 
=\epsilon^{-1}T_{\ov X/\ov S}(\log D)\subset \sA_E;\\
\notag\sA_{\lambda_E}(\log W) =
\epsilon^{-1}T_{\ov S}(\log W)\subset \sA_{\lambda_E},
\end{gather}
where ``$\epsilon$" denotes symbol maps.
\end{notation}
The $\,^{tr}\sA^{-1}_E$ in Section 3 admits a generalization
$\,^{tr}\sA^{-1}_E(\log D)$ (cf. \cite{TT}, p. 593) such that
\begin{gather}
0\to \omega_{\ov X/\ov S}\to \,^{tr}\sA^{-1}_E(\log D)\to 
\sA_{E/\ov S}(\log D)\to 0,
\end{gather}
is exact.  Furthermore with $\,^{tr}\sA^{\bullet}_E$ 
in Section 3 replaced by $\,^{tr}\sA^{\bullet}_E(\log D)$
(with $\sA_{E,\pi}(\log D)$ in degree 0), one has

\begin{prop} \label{prop6.3} {\rm (cf. \cite{TT})}
There is a canonical isomorphism
$$R^0\pi_*^{tr}\sA^{\bullet}_E(\log D)\cong
\sA_{\lambda_E}(\log W)$$ that extends
Proposition \ref{prop4.1}.
\end{prop}

Now we come back to moduli situation. 
Recall Notation \ref{prop2} and the Kodaira-Spencer map (cf. \cite{TUY}, Remark 3.2.7)
\begin{gather}\label{eq6.2}
{\rm KS}:{\ov f}^*T_{\ov\sM}(B)\to R^1\pi_*T_{\ov X/\ov S}(\log D).
\end{gather}
As the same as the situation of smooth curves, we have 
\begin{lem}\label{prop9.9} The morphism
$\ov f: \ov S\to\ov\sM$ is smooth and 
$$T_{\ov S/\ov\sM}=R^1\pi_*\sE nd^0(E).$$
\end{lem}

\begin{proof}
When $C$ is irreducible, its fibre 
${\ov f}^{-1}([C])$ is the moduli space of stable bundles with fixed
determinant $\sN|_C$. When $C$ is reducible, ${\ov f}^{-1}([C])$
have a few disjoint irreducible components and
each component consists of bundles with a fixed determinant that
coincides with $\sN|_C$ outside the node of $C$ (cf. \cite{S}).
\end{proof} 

\begin{prop} \label{prop6.4}
\item i) Assume
the Kodaira-Spencer map ${\rm KS}$ \eqref{eq6.2} is injective.
Then $$R^0\pi_*(\sE nd(E)^{-1}\to \sA_{E, \pi}(\log D))
\to R^0\pi_*(\sE nd(E)\to \sA_{E,\pi}(\log D))\to 0$$
is canonically isomorphic to
$$\sD^{\le 1}_{\ov S/\ov\sM}(\lambda_E)\to T_{\ov S/\ov\sM}\to 0.$$
\item ii) Assume that the ${\rm KS}$ \eqref{eq6.2}
is an isomorphism. 
Then $$ R^1\pi_*\sA^0_{E/\ov S}(\log D)\cong T_{\ov S}(\log W)$$
canonically, where $\sA^0_{E/\ov S}(\log D)= \sA_{E/\ov S}(\log D)/ \sO_{\ov X}$.
\end{prop}

The following is left to the reader.

\begin{lem} \label{prop6.7}
\begin{gather}\notag
(\sA_{E/\ov S}(\log D)\to \sA_{E,\pi}(\log D))\cong_{\rm q.i.}
(T_{\ov X/\ov S}(\log D)\to T_{\pi}(\log D))
\\ \notag \qquad
\cong_{\rm q.i.}\pi^{-1}T_{\ov S}(\log W)
\end{gather} as quasi-isomorphisms.
\end{lem}

We prove now Proposition \ref{prop6.4}.
\begin{proof}
i) From the exact sequence
\begin{gather}
\notag
0\to
(\sE nd(E)^{-1}\to \sA_{E,\pi}(\log D))
\to\, ^{tr}\sA^{\bullet}_E(\log D)\to T_{\ov X/\ov S}(\log D)[1]\to 0,
\end{gather}
and passage to cohomology
\begin{gather}
\notag 0\to
R^0\pi_*(\sE nd(E)^{-1}\to \sA_{E,\pi}(\log D))
\to \sA_{\lambda_E}(\log W)\xrightarrow{\epsilon}
\\ \notag\qquad R^0\pi_*T_{\ov X/\ov S}(\log D)[1]\to \cdots,
\end{gather}
one has, via the injectivity of ${\rm KS}$,
that $\epsilon^{-1}(0)$ has symbolic part in $T_{\ov S/\ov\sM}$.
This gives one of the isomorphisms in i)
(the one with $D^{\le 1}_{\ov S/\ov\sM}(\lambda_E)$).
Further, by \begin{gather}
\notag 0\to\omega_{\ov X/\ov S}[1]\to
(\sE nd(E)^{-1}\to \sA_{E,\pi}(\log D))
\to (\sE nd(E)\to \sA_{E, \pi}(\log D))\to 0
\end{gather}
and $R^0\pi_*\omega_{\ov X/\ov S}[1]\cong \sO_{\ov S}$ it follows $$
R^0\pi_*(\sE nd(E)^{-1}\to \sA_{E,\pi}(\log D))
\to R^0\pi_*(\sE nd(E)\to \sA_{E,\pi}(\log D))$$
is nothing but the symbol map, completing
the asserted isomorphisms.
ii)  Write $B^{\bullet}(\log D)$
for $(\sE nd(E)\to \sA_{E,\pi}(\log D))$
and $\sA^{\bullet}_{E,\pi} (\log D)$
for L.H.S. of Lemma \ref{prop6.7}.
The following exact sequence
\begin{gather}\label{eq7.6.1} 0\to B^{\bullet}(\log D)\to
\sA^{\bullet}_{E,\pi}(\log D)\to T_{\ov X/\ov S}(\log D)[1]\to 0
\end{gather}
projects to
\begin{gather}\label{eq7.6.2}
0\to \sE nd^0(E)[1]\to \sA^0_{E/\ov S}(\log D)[1]\to
T_{\ov X/\ov S}(\log D)[1]\to 0
\end{gather}
Computing direct images of \eqref{eq7.6.1} and \eqref{eq7.6.2},
by i) just proved
and Lemma \ref{prop6.7}, yields that
$R^0\pi_*$ of \eqref{eq7.6.2}
should be isomorphic to
$$0\to T_{\ov S/\ov\sM}\to T_{\ov S}(\log W)\to {\ov f}^*T_{\ov\sM}(B)\to 0,$$
implying the assertion.
\end{proof}

\begin{remark}\label{prop6.4.r}
i) $R^1\pi_*(\sG_E)\cong
R^0\pi_*(\sE nd(\E)^{-1}\to \sA_{\E,\pi}(\log D))
\cong
\sD^{\le 1}_{\ov S/\ov\sM}(\lambda_{\E})$ holds, cf. Theorem \ref{prop5.4}
and Remark \ref{prop5.5}.
ii) For the family $\ov\sC\to \ov\sM$
there is a Kodaira-Spencer map
(cf. \cite{TUY}, 3.1-2)
\begin{gather}\label{eq6.2.1}
\rho_b:T_{\ov\sM,b}\to {\rm Ext}^1_{\sO_{\ov\sC_b}}(\Omega_{\ov\sC_b},
\sO_{\ov\sC_b}),
\,\, b\in \ov\sM
\end{gather}
The family $\ov\sC\to \ov\sM$ is a {\it local universal family} if $\rho_b$ is
an isomorphism at each $b\in \ov\sM$.  If $\ov\sC\to \ov\sM$ is a local
universal family (cf. \cite{TUY}, Theorem 3.1.5 for
the existence of such a family),
then ${\rm KS}$ \eqref{eq6.2} is an isomorphism
(cf. \cite{TUY}, Theorem 3.2.6).
\end{remark}

Combining the above with 
the 2nd half of Section 2, 
we are now ready to generalize 
Theorem \ref{prop5.4} in the context 
of {\it log geometry}.

\begin{thm}\label{prop6.6}
Suppose ${\rm KS}$ in \eqref{eq6.2} is an isomorphism 
({\rm cf}. ii) of Remark \ref{prop6.4.r}). Then, over any
affine open set $U_i$ of $\ov\sM$, there is an  
 isomorphism 
$$\phi_i:R^1\pi_*S(\sG_E) \cong 
\sD_{\ov S}^{\le 1}(\lambda_{\E})(\log W).$$ 
\end{thm}
\begin{proof} 
Note the above assumption for   
${\rm KS}$ is, 
via ii) of Proposition \ref{prop6.4},
in correspondence to ${\rm KS}_S$ 
in Theorem \ref{prop5.4}. 
It follows that all key ingredients in proof of Theorem
\ref{prop5.4} admit corresponding counterparts 
for log geometry, such as Proposition \ref{prop6.3} 
and Proposition \ref{prop6.4}. 
Thus the generalization of Theorem \ref{prop5.4}
in {\it log} context is immediate.  
\end{proof}

To complete this paper, we prove the coherence of $f_{T*}\Theta^k$, for
which our proof need a result on the density of
locally free sheaves.

\begin{prop}\label{prop9.1} The fibre of
$f_Z: Z\to\ov\sM$ at any point of $B=\ov\sM\setminus\sM$ has a dense
open set of locally free sheaves.
\end{prop}

\begin{proof} Let $X_0$ be a fibre of $\sC\to\ov\sM$ at $0\in B$.
If $X_0$ is reducible,
the lemma is known (see \cite{S}). Thus we assume that $X_0$ is irreducible.
Let $\sU_0$ be the moduli
space of semistable sheaves of rank $r$ and degree $d$ (without fixed
determinant) on $X_0$. We need to show that
$Z_0:=f_Z^{-1}(0)$ ($\subset \sU$)
contains a dense set of locally free sheaves with the fixed
determinant $\sN_0$. Let $J(X_0)$ be the Jacobian of $X_0$, which consists
of line bundles of degree $0$ (thus non-compact for the singular curve
$X_0$). Then we have a morphism
$$\phi_0: Z_0\times J(X_0)\to \sU_0,\quad \text{where $\phi_0(E,L)=E\otimes 
L$}.$$
Now we prove that $\phi_0$ has fibre dimension at most $1$, namely,
for any $[F_0]\in \sU_0$, the fibre
$$\phi_0^{-1}([F_0])=\{(F,L)|F\otimes L=F_0\}\subset Z_0\times J(X_0)$$
has at most dimension $1$ (for simplicity, we assume that
$X_0$ has only one node $x_0$).
One can check that for any $[F]\in Z_0$ it satisfies 
$$\frac{\wedge^rF}{\rm torsion}\subset \sN_0\quad \text{and}\quad
m_{x_0}^{r(F)}\sN_0\subset\frac{\wedge^rF}{\rm torsion}$$
where $m_{x_0}$ is the ideal sheaf of the node
$x_0\in X_0$ and $r(F)=d-{\rm deg}(\frac{\wedge^rF}{\rm torsion})$.
Let $\rho: \widetilde{X_0}\to X_0$ be the normalization and
$\rho^{-1}(x_0)=\{x_1,x_2\}$. Then the above condition implies that
$$\frac{\rho^*(\wedge^rF)}{\rm torsion}=\rho^*\sN_0(-n_1x_1-n_2x_2),\quad
n_1+n_2\le 2r(F),\,n_1\ge 0,n_2\ge 0,$$
where $r(F)=d-{\rm deg}(\frac{\wedge^rF_0}{\rm torsion})=r(F_0)$
since $F\otimes L=F_0$ and thus
$$\frac{\wedge^rF_0}{\rm torsion}=\frac{\wedge^rF}
{\rm torsion}\otimes L^r.$$
Thus, for any $(F,L)\in\phi_0^{-1}(F_0)$, $L$ has to satisfy
$$(\rho^*L)^r=(\rho^*\sN_0)^{-1}(n_1x_1+n_2x_2)\otimes
\frac{\rho^*\wedge^rF_0}{\rm torsion}$$
which is a finite set since there are only finitely many choices of $(n_1,n_2)$.
The pullback map $\rho^*:J(X_0)\to J(\widetilde {X_0})$ has
$1$-dimensional kernel. Thus $\dim(\phi_0^{-1}([F_0]))\le 1.$
We shall now prove the density of locally free sheaves.

Let $Z^i_0$ be an irreducible component of $Z_0$, then
\begin{gather}
\label{eqdim}
\dim(Z_0^i\times J(X_0))\ge\dim(Z_{\eta}\times J(X_0))=\dim(\sU_0). 
\end{gather}
If $Z_0^i$ contains no locally
free sheaf, then $\phi(Z_0^i\times J(X_0))$ falls into the
subvariety $\sU_0^n\subset\sU_0$ of non-locally free sheaves.
$\sU^n_0$ has a dense open set $\sU^n_0(1)$ consisting of 
torsion free sheaves
$F$ of the following type, said to be type $1$.  Namely,
$$F\otimes\hat\sO_{X_0,x_0}=\hat\sO_{X_0,x_0}^{(r-1)}\oplus \hat m_{x_0}.$$
If $\phi((F,L))=F_0\in \sU^n_0(1)$, then $F$ is also of type $1$
(tensoring a line bundle do not change its type). By Remark 8.1 of
\cite{NS},
$$L_0=\frac{\wedge^rF}{\rm torsion}$$
is a torsion free (but non-locally free) sheaf of degree $d$. But 
$L_0\subset \sN_0$, thus $L_0=\sN_0$ since
they have the same degree, a contradiction with that $L_0$ is 
not locally free.  
\end{proof}

\begin{thm}\label{prop9.2}  i) $f_{T*}\Theta^k$ is coherent;
ii) $f_{T*}\Theta^k={\ov f}_*\Theta^k$ if either $g\ge 3$ or $r\ge 3$.
\end{thm}

\begin{proof}
Let $\iota:\ov Z\to Z$ be the normalisation of $Z$ and
$\ov\Theta=\iota^*(\Theta)$. Write
$f_{\ov Z}:\ov Z\to\ov\sM$ and $f_{\iota^{-1}(T)}:\iota^{-1}(T)\to\ov\sM$.
Then $\iota^{-1}(T)\cong T$ since
$T$ is normal, and
$\sF:=f_{T*}\Theta^k\cong f_{\iota^{-1}(T)*}({\ov\Theta}^k).$
On the other hand, since each fibre of $f_Z:Z\to\ov\sM$ contains
a dense open set of locally free sheaves, we have
${\rm codim}(Z\setminus T)\ge 2$. Thus
$\ov Z\setminus \iota^{-1}(T)=\iota^{-1}(Z\setminus T)$
has codimension at least $2$ since $\iota:\ov Z\to Z$ is a finite map.
By Hartogs type extension theorem,
$$\sF\cong f_{\iota^{-1}(T)*}({\ov\Theta}^k)\cong f_{{\ov 
Z}*}({\ov\Theta}^k),$$
which is coherent, hence i). The claim ii) that $f_{T*}\Theta^k={\ov 
f}_*\Theta^k$
follows also from the Hartogs type theorem because
$T$ is normal and $T\setminus \ov S$ is of
codimension at least $2$ when $g>2$ or $r>2$ (cf. \cite{H}).
\end{proof}

Finally, we prove a lemma, which is not needed for
this paper. But we expect that it will be useful in the future
since it gives the relationships of $S(\sG_E)$, $\sG_{E^*}$ and
$\sA_{E^*/S}$. We remark that the morphism $\ov{\rm res}$ in   
Proposition \ref{props6.12} is induced by pairings in the lemma.
 
\begin{lem}\label{prop5.11}
There are canonical isomorphisms:
\item i) $S(\sG_E)/\omega_{X/S}
\cong (\sG_{E^*})^*$,
\item ii) $(\sG_E)^*\cong \sA^0_{E^*/S}$
where $\sA^0_{E^*/S}=\sA_{E^*/S}/\sO_X$.
\end{lem}

\begin{proof}
i) One constructs a non-degenerate pairing
$$\sG_{E^*}\times
S(\sG_E)/\omega_{X/S}
\to \sO_X$$
using local description \eqref{eq5f1} (also cf. \cite{BS}). We define 
\begin{gather}\label{eqp}
<s_1,s_2>:= \chi\nu+{1\over 2}\sum_a
(\nu_a+\mu_a){\rm Trace}(J_a\cdot\, ^tB)
\end{gather}
for $s_1=(0,B,\nu)$ and
$$s_2=\left(\aligned&\chi\sum_a(0,J_a, 0)\otimes (0,J_a,0)+
\sum_a\mu_a(0,J_a,0)\otimes(0,0,1)\\&+\sum_a\nu_a(0,0,1)\otimes(0,J_a,0)
+w(0,0,1)\otimes(0,0,1)\endaligned\right)\otimes\partial_t.$$
One sees that
$\omega_{X/S}$ is contained
in the kernel of the pairing \eqref{eqp}. The pairing is obviously
invariant under coordinates change $y(t)$. If the trivialization of $E$
is changed by a gauge $g$, then $E^*$ is changed by a gauge $(^tg)^{-1}$.  
Thus the verification of \eqref{eqp} being $g$-invariant 
is easily reduced to an identity 
\begin{gather}\label{eq5f8}    
\sum_a{\rm Tr}(-J_ag^{-1}g')\cdot{\rm Tr}(J_a\cdot{^tB}) 
={\rm Tr}(B(^tg)((^tg)^{-1})'). 
\end{gather} 
The L.H.S. of \eqref{eq5f8} equals ${\rm Tr}(-g^{-1}g'\cdot {^tB})$ 
($B$ being traceless).       
The R.H.S. of \eqref{eq5f8}, after transposition, is 
$${\rm Tr}((g^{-1})'g\cdot{^tB})=  
{\rm Tr}\big((-g^{-1}g'g^{-1})g\cdot ^tB\big).$$
  
ii) Define a non-degenerate pairing
\begin{gather}\label{eqp2}
\sG_E\times\sA^0_{E^*/S}\to \sO_X
\end{gather}
by $<(0,B_1,\nu), (\chi, B_2)>=
\nu\chi+{\rm Trace}(B_1\cdot ^{ t}B_2)$. We check that
it is independent of choices of coordinates and gauges.
\begin{gather}\label{eq5f9}   
(0, B_1, \nu):=s_1
\buildrel{g}\over \longrightarrow 
\big(0, gB_1g^{-1}, {\rm Tr}(-B_1g^{-1}g')+\nu\big); 
\\
(\chi, B_2):=s_2 \buildrel{^tg^{-1}}\over \longrightarrow 
\Big(\chi, -\chi\big((^tg^{-1})'\cdot{^tg}\big)+(^tg^{-1})B_2\cdot^tg\Big);  
\notag
\\
s_1\buildrel{y(t)}\over \longrightarrow 
(0, B_1, \nu y'); \qquad s_2    
\buildrel{y(t)}\over \longrightarrow 
\big(\chi {y'}^{-1}, B_2\big).   
\notag
\end{gather}  
The $y(t)$-change is obvious.  
For $g$-change,  we have 
\begin{gather}\label{eq5f10}   
\chi \nu + {\rm Tr}(B_1\cdot ^tB_2)
\longrightarrow 
\chi \nu + \chi {\rm Tr}(-B_1g^{-1}g')+ 
\\  
{\rm Tr}\Big(gB_1g^{-1}\cdot ^t\big(-\chi ((^tg^{-1})'
\cdot {^tg}\big)\Big)
+{\rm Tr}\Big(gB_1g^{-1}\cdot ^t\big(^tg^{-1}B_2\cdot ^tg\big)\Big).
\notag
\end{gather}    
The 1st (resp. 2nd) term in the last line equals 
\begin{gather}\label{eq5f11}
-\chi{\rm Tr}\big(gB_1g^{-1}\cdot g(g^{-1})'\big)
=-\chi{\rm Tr}\big(gB_1g^{-1}\cdot g(-g^{-1}g'g^{-1})\big)
\\ 
=\chi{\rm Tr}\big(g(B_1g^{-1}g')g^{-1}\big)
=\chi{\rm Tr}(B_1g^{-1}g')  
\notag
\\
\big({\rm resp.}\quad  {\rm Tr}(gB_1\cdot ^tB_2g^{-1})=
{\rm Tr}(B_1\cdot ^tB_2)\big).   
\notag
\end{gather}  
It follows from \eqref{eq5f10} and 
\eqref{eq5f11} that the pairing $\chi\nu+\Tr(B_1\, ^tB_2)$ 
in \eqref{eqp2} is globally defined.      
\end{proof}

\bibliographystyle{plain}
\renewcommand\refname{References}

\end{document}